\documentclass{amsart}
\usepackage{amsmath}
\usepackage{accents}
\usepackage{paralist}
\usepackage[misc]{ifsym}
\usepackage{epsfig} 
\usepackage{epstopdf} 
\usepackage{mathrsfs}
\usepackage{txfonts}
\usepackage{amsfonts, amsthm, amssymb}
\usepackage{mathtools} 
\usepackage{graphicx} 
\usepackage{enumerate}
\usepackage{enumitem}
\usepackage{color}
\usepackage{fancyhdr}
\usepackage{diagbox}
\usepackage{threeparttable}
\usepackage{tabularx}
\usepackage{booktabs}
\usepackage{comment}
\usepackage{tikz}
\usepackage{hyperref}       
\usepackage{url}           

\allowdisplaybreaks

\numberwithin{equation}{section}
\newtheorem{theorem}{Theorem}[section]
\newtheorem{corollary}[theorem]{Corollary}

\newtheorem{lemma}[theorem]{Lemma}
\newtheorem{proposition}[theorem]{Proposition}

\theoremstyle{definition}
\newtheorem{definition}[theorem]{Definition}
\newtheorem{remark}[theorem]{Remark}

\newcommand{\silentfootnote}[1]{
	\begingroup
	\renewcommand{\thefootnote}{}
	\def\@makefnmark{}
	\footnote{#1}
	\addtocounter{footnote}{-1}
	\endgroup
}

\usetikzlibrary{calc}
\allowdisplaybreaks

\selectfont

\title[Regularity and Symmetry of Solutions to Mixed  Parabolic Equations]{Regularity and Asymptotic Radial Symmetry of Solutions to Mixed Local-Nonlocal Parabolic Equations}

\author[Fuwei Cheng]{Fuwei Cheng}
\address{School of Mathematical Sciences, Beijing Normal University, No. 19, XinJieKouWai St., HaiDian District, Beijing 100875, P. R. China}
\email{fwcheng@mail.bnu.edu.cn}

\author[Xifeng Su]{Xifeng Su}
\address{School of Mathematical Sciences, Laboratory of Mathematics and Complex Systems (Ministry of Education)\\
	Beijing Normal University,
	No. 19, XinJieKouWai St., HaiDian District, Beijing 100875, P. R. China}
\email{xfsu@bnu.edu.cn, billy3492@gmail.com}

\begin{document}
\begin{abstract}
		We study the regularity and the long-time asymptotics of the solution to the following mixed local-nonlocal parabolic problem
		\[
		\left\{
		\begin{alignedat}{2}
			\partial_t u - \Delta u + (-\Delta)^s u &= f \qquad &&\text{in }\Omega\times(0, T),\\
			u&=0 \qquad &&\text{in }\bigl(\mathbb{R}^{n}\setminus\Omega\bigr)\times[0, T),\\
			u(\,\cdot\,,0)&=u_{0} \qquad &&\text{in }\Omega,
		\end{alignedat}
		\right.
		\]
		where $\Omega\subset\mathbb R^n$ is a bounded domain, $s\in(0,1)$ and $(-\Delta)^s$ denotes the fractional Laplacian.
		
		We first develop the well-posedness theory for the problem on finite time intervals ($0<T<+\infty$). 
		Existence and uniqueness of weak solutions are established in the natural energy space, followed by interior $H^{2,1}$-regularity and global $L^\infty$-bounds.
		
		Then, under the restriction $s\in(0,1/2)$ and $0<\alpha\le 1-2s$, we prove that weak solutions become classical ones, obtaining $C^{2+\alpha,1+\alpha/2}$-regularity up to the boundary via the method of continuity and a perturbation argument.
		The accompanying a priori estimate guarantees that these classical solutions can be extended globally in time, thereby allowing us to consider the parabolic problem on the whole half-line ($T=+\infty$).
		
		Finally, for the special case $\Omega=B_1(0)$ and $f=0$, we employ the moving-plane method together with the previously established boundary regularity to show that every positive classical solution is attracted, as $t\to +\infty$, by a steady state that is radially symmetric and strictly decreasing.
\end{abstract}

\thanks{X. Su is supported by the National Natural Science Foundation of China (Grant No. 12371186, 12571207).}	

\maketitle
\noindent{\it Keywords:} mixed local-nonlocal parabolic equation; second-order H\"older regularity; regularity up to the boundary; method of continuity; moving plane method; perturbation method; asymptotic radial symmetry.

	
\section{Introduction}
We focus first on the study of the regularity results 
 of solutions to the mixed local-nonlocal parabolic equation
\begin{equation}\label{problem}
	\partial_t u - \Delta u + (-\Delta)^s u = f \quad \text{in } Q_T:=\Omega\times(0,T),
\end{equation}
subject to the exterior Dirichlet condition
\begin{equation}\label{bdry-ext}
	u=0 \quad \text{in }\bigl(\mathbb{R}^{n}\setminus\Omega\bigr)\times[0,T),
\end{equation}
and the initial condition
\begin{equation}\label{bdry-initial}
	u(\,\cdot\,,0)=u_{0} \quad \text{in }\Omega.
\end{equation}
Throughout the paper $\Omega\subset\mathbb{R}^{n}$ is a bounded domain, $s\in(0,1)$ and $(-\Delta)^{s}$ stands for the fractional Laplacian defined by the singular integral
\[
(-\Delta)^{s}u(x,t)=c_{n,s}\,\mathrm{P.V.}\int_{\mathbb{R}^{n}}\frac{u(x,t)-u(y,t)}{|x-y|^{n+2s}}\,dy,
\]
with $c_{n,s}>0$ a universal normalization constant.

The mixed differential-pseudo-differential operator  
\[
\mathcal{L}:=-\Delta+(-\Delta)^{s}, \qquad s\in(0,1),
\]  
is the generator of the superposition of a Brownian motion and a $2s$-stable L\'evy process, and governs models for instance in population dynamics \cite{DV21,DPV23,PV18}, plasma physics \cite{BD13}, as well as finance and stochastic control \cite{MP96}.

Elliptic theory involving $\mathcal{L}$ is now well-developed, see e.g.   existence/non-existence results \cite{SVWZ24,RS15}, regularity analysis \cite{BDVa22,SVWZ25}, Harnack inequalities \cite{F09,CKSV12,GK22}, viscosity solutions theory \cite{EK05,GC08}, and references therein.

The counterpart of parabolic equations $\partial_{t}u+\mathcal{L}u=0$ have been investigated more recently. A weak Harnack inequality was proved in \cite{GK23}; quasilinear extensions with $p$-Laplacian type diffusion were studied in \cite{FSZ22,GK24}.

In this paper, we will  deal with the mixed parabolic problem \eqref{problem}-\eqref{bdry-initial}. 
Let us first develop the weak solution theory.
Note that the candidate space of weak solutions of \eqref{problem} would be the parabolic energy space\footnote{See Definition~\ref{parabolicSobolev}.}
\begin{align*}
	X^{1,1}(R_{T}) = \left \{ u:
	\|u\|_{X^{1,1}(R_{T})} := \|u\|_{L^{2}(R_{T})} + \|D_{t}u\|_{L^{2}(R_{T})} + \|\nabla u\|_{L^{2}(R_{T})} + \int_{0}^{T} [u]_{s}^{2}\, dt<\infty  \right\},
\end{align*}
with \(R_T=\mathbb{R}^{n}\times(0,T)\).
By the continuous embedding $H^{1}(\mathbb{R}^{n})\hookrightarrow H^{s}(\mathbb{R}^{n})$ (see e.g. \cite[Proposition 2.2]{DNPV12}), we obtain the continuous embedding $X^{1,1}(R_T)\hookrightarrow H^{1,1}(R_T)$. 

To incorporate the homogeneous exterior Dirichlet condition \eqref{bdry-ext}, we thus consider the space \( X^{1,1}_0(Q_T)\) defined as the closure of smooth functions vanishing on \((\mathbb{R}^n\setminus\Omega)\times(0,T)\) in the norm of \(X^{1,1}(R_T)\).

Leveraging the elliptic theory of \cite{BDVa22}, we obtain the existence and uniqueness results by Galerkin approximation, derive the interior Sobolev regularity via truncation techniques and energy estimates, and establish the $L^{\infty}$-boundedness through De Giorgi iteration.
	
\begin{theorem}\label{Sobolevregularity}
	Let $\partial\Omega\in C^1$. Assume $f\in L^{2}(Q_{T})$ and $u_0 \in L^2(\Omega)$, and then problem \eqref{problem}--\eqref{bdry-initial} admits a unique weak solution $u\in X^{1,1}_0(Q_{T})$. 
	Moreover,
	\begin{itemize}			
	 \item [(i)]  	 
	if, in addition, $u_0\in H^1(\Omega)$, then $u\in H^{2,1}(Q_{T}^{\delta})$ for any $\delta>0$,
	where  the interior cylinder  
	$
	Q_{T}^{\delta}:=\bigl\{x\in\Omega:\operatorname{dist}(x,\partial\Omega)>\delta\bigr\}\times(0,T),
	$ 
	 and  
	\[
	\|u\|_{H^{2,1}(Q_{T}^{\delta})}\le C\bigl(\|u\|_{H^{1,0}(Q_{T})}+\|f\|_{L^{2}(Q_{T})}+\|u_{0}\|_{H^1(\Omega)}\bigr),
	\]  
	where $C=C(n,\delta,T,\Omega)$ is a positive constant.
	\smallskip
	
	\item [(ii)] if $f\in L^{\infty}(Q_{T})$ and $u_0 \in L^\infty(\Omega)$, then $u\in L^{\infty}(Q_{T})$ and
	\[
	\|u\|_{L^{\infty}(Q_{T})}\leq C\bigl(\|f\|_{L^{\infty}(Q_{T})}+\|u_{0}\|_{L^{\infty}(\Omega)}\bigr)
	\]
	with some positive constant $C=C(n,T,\Omega)$.
	\end{itemize}
\end{theorem}


Although interior Sobolev regularity of the weak solution of problem \eqref{problem}--\eqref{bdry-initial} is now available, we deliberately pursue its H\"older boundary regularity rather than the  Sobolev boundary regularity. 
The reason is twofold. 
\begin{itemize}
\item On the one side, it is rather a challenge to show directly the Sobolev boundary regularity in the present setting since several methods of proof fail  due to the fact that 
the standard boundary-flattening conflicts with the exterior condition required by $(-\Delta)^s$. For instance, the upgraded global diffeomorphism approach of \cite{DFM24} does not work easily because of substantial complexity.

\item On the other side, once the H\"older boundary regularity\footnote{We actually establish $C^{2+\alpha,1+\alpha/2}$-regularity for solutions of problem~\eqref{problem}--\eqref{bdry-initial} in Theorem~\ref{thm: reg up to bound} below.} is obtained, as an application, we are able to show the asymptotic behavior of the solutions to the parabolic problem (see Theorem~\ref{thm:symmetry}) by the 
the moving-plane argument, which demands pointwise control of second-order derivatives and continuity at the boundary.
\end{itemize}

Our next goal is to establish the second order H\"older regularity up to the boundary for solutions to the mixed local-nonlocal parabolic problem~\eqref{problem}-\eqref{bdry-initial}.
We would point out that Bin Shang and Chao Zhang \cite{BZ22} proved the local H\"older continuity for the mixed $p$-Laplacian evolution  
\[
\partial_t u-\Delta_p u+(-\Delta_p)^s u=0, \quad \text{in } \Omega \times (0,T)
\]  
via expansion of positivity and De Giorgi iteration.  
Stuti Das \cite{Das24} treated the linear case $p=2$ and derived interior gradient H\"older regularity by perturbation against the classical heat equation.  

We now advance the regularity theory by establishing second-order H\"older continuity up to the boundary for the linear case $p=2$.
In fact, note that in the case $s\in (0, 1/2)$, by treating the nonlocal term as a lower-order perturbation of the classical heat operator, 
we apply the method of continuity together with Schauder theory to achieve full boundary regularity.

\begin{theorem}[$C^{2+\alpha,1+\alpha/2}$-regularity up to the boundary]\label{thm: reg up to bound}
	Let 
	{\(s\in(0,1/2)\) and \(0<\alpha\leq 1-2s\),}
	$\Omega\subset\mathbb R^n$ be a bounded $C^{2+\alpha}$-domain, $f\in C^{\alpha,\alpha/2}(\overline Q_T)$ and $u_0\in C^{2+\alpha}(\overline\Omega)$.  
	If $u$ is a weak solution of problem \eqref{problem}-\eqref{bdry-initial}, then  
	\(
	u\in C^{2+\alpha,1+\alpha/2}(\overline Q_T)
	\) 
	and   
	\begin{equation}\label{reg-up-to-bd}
		\|u\|_{C^{2+\alpha,1+\alpha/2}(\overline Q_T)}
		\le C\left(\|f\|_{C^{\alpha,\alpha/2}(\overline Q_T)}+\|u_0\|_{C^{2+\alpha}(\overline\Omega)}+\|u\|_{L^\infty(Q_T)}\right),
	\end{equation}
	where $C>0$ is a constant depending only on $n,s,\alpha,\Omega,T$.
\end{theorem}

We point out that the proof for Theorem~\ref{thm: reg up to bound} is self-contained and does not rely on the interior second-order theory developed in the Sobolev framework as in Theorem~\ref{Sobolevregularity}.
We remark that the restriction $s<1/2$ is somewhat unavoidable in view of the stationary case \cite[Remark 4.6]{CSZ25}.

As a consequence, one can obtain the following global-in-time existence of classical solutions for smooth data.
\begin{corollary}\label{cor:global-exist-reg}
	Assume that all the conditions in Theorem~\ref{thm: reg up to bound} are satisfied.
	Let $f\in C^{\alpha,\alpha/2}(\overline Q_\infty)$.
	Then the problem
	\begin{equation}\label{eq:globalintime}
	\left\{
	\begin{alignedat}{2}
		\partial_t u-\Delta u+(-\Delta)^s u &=f &\quad&\text{in }\Omega\times(0, +\infty),\\
		u&=0 &&\text{in }(\mathbb{R}^n\setminus\Omega)\times[0, +\infty),\\
		u(\cdot\,,0)&=u_0 &&\text{in }\Omega
	\end{alignedat}
	\right.
	\end{equation}
	possesses a solution $u\in C^{2+\alpha,1+\alpha/2}(\overline Q_\infty)$.
\end{corollary}
	
This global well-posedness in turn motivates the study of long-time behavior of the classical solutions to \eqref{eq:globalintime}, such as the asymptotic radial symmetry and monotonicity.
 
Note that the study of symmetry in parabolic equations has a long history, beginning with the seminal work of Gidas-Ni-Nirenberg \cite{GNN79} on the classical heat equation, which inspired extensive subsequent literature such as \cite{HP94,Bab94,Fol11,SnW12}.  
Parallel results for purely nonlocal parabolic operators were later obtained in e.g. \cite{JW14,CWNH21,WC22}. 
 
We aim to extend the above conclusions to the present mixed setting. 
More concretely, following the dynamical perspective, 
we study the $\omega$-limit set 
\[
\omega (u)= \{ \varphi = \lim u(\cdot,t_k) \text{ for some } t_k \rightarrow +\infty\}
\]
which collects all possible limiting profiles. 

 In fact, one can exploit the Hille-Yosida theorem (see \cite[Thm.~4.9]{Gri09}) to ensure that the operator $-\mathcal{L}$ generates a strongly continuous (analytic) semigroup
\(\bigl(e^{-t\mathcal{L}}\bigr)_{t\ge 0}\) on $L^{2}(\mathbb{R}^{n})$.  	
Writing $\mathcal{T}(t,u_{0}):=e^{-t\mathcal{L}}u_{0}$, we obtain a dynamical system
\[\bigl\{\mathcal{T}(t,u_{0}):t\ge 0\bigr\}\]
whose trajectories $t\mapsto u(\cdot,t)$ are precisely the mild solutions of the parabolic equation.

Owing to the boundary \(C^{2+\alpha,1+\alpha/2}\)-regularity established in Theorem \ref{thm: reg up to bound}, the solution \(u\) is equicontinuous whenever \(f\in C^{\alpha,\alpha/2}(\overline{Q}_{T})\) and \(u_{0}\in C^{2+\alpha}(\overline{\Omega})\).  
Consequently, \cite[Proposition 53.3]{QS19} ensures that the \(\omega\)-limit set \(\omega(u)\subset C^{0}(\mathbb{R}^{n})\) is nonempty, compact and attracts the trajectory in the uniform norm:
\[
\operatorname{dist}_{C^{0}(\mathbb{R}^{n})}\, \bigl(u(\cdot,t),\omega(u) \bigr)\to 0
\qquad\text{as } t\to +\infty.
\]


Therefore, we obtain the following asymptotic behavior for the mixed local-nonlocal parabolic equations.
\begin{theorem}[Asymptotic radial symmetry and monotonicity]\label{thm:symmetry}
	Let $s\in(0,1/2)$, \(\Omega=B_1(0)\) and \(f=0\).  
	Assume $u \in  C^{2+\alpha,1+\alpha/2}(\overline Q_\infty)$ is a positive solution to   
	\begin{equation}\label{symmetry problem}
	\left\{
	\begin{aligned}
		\partial_t u-\Delta u+(-\Delta)^s u&=0 &&\text{in } B_1(0)\times(0, +\infty),\\
		u&=0 &&\text{in } \left(\mathbb{R}^n \setminus B_1(0) \right) \times[0, +\infty).
	\end{aligned}
	\right.
	\end{equation} 
	Then, any steady state $\varphi\in\omega(u)$ is either identically zero or strictly radially decreasing. 
\end{theorem}

\begin{remark}
We point out that, for fixed initial data \(u(\cdot, 0) = u_0 > 0\), the problem \eqref{symmetry problem} admits a unique positive solution that is asymptotically radially symmetric, 
which follows from Corollary~\ref{cor:global-exist-reg} and the weak maximum principle (see Theorem~\ref{thm:max_principle}).
\end{remark}


For simplicity and in order to get right to the point, we will only consider the specific case when $f=0$ which may highlight the genuine geometric effect of $\partial_t+\mathcal{L}$.
Extensions to general nonlinearities $f(t, x, u)\not\equiv 0$, such as $u^p-u$ with $p>1$, will be pursued in forthcoming work. 

The paper is structured as follows.  
Section~\ref{sec: Preliminaries} collects the notations and function spaces that will be used throughout, and gives a basic weak maximum principle adapted to the mixed-order operator. 
In Section~\ref{sec: Regularity of weak solutions}, we develop the well-posedness and regularity theory for the mixed parabolic problem:  
we begin by establishing the existence and uniqueness of the weak solution; 
then, our analysis proceeds from basic Sobolev regularity to fine pointwise estimates -- interior $H^{2,1}$-regularity, global $L^\infty$-bounds and the higher-order $C^{2+\alpha,1+\alpha/2}$ boundary regularity;  
finally, we obtain the global-in-time existence of the classical solution.
Based on the above results, Section~\ref{sec: Asymptotic Radial Symmetry of Classical Solutions} investigates the long-time behavior of the positive classical solution to the homogeneous equation on the unit ball.  

\section{Preliminaries}\label{sec: Preliminaries}
\subsection{Useful notations} 
We introduce some standard notations that will be used throughout the paper. 
\begin{itemize}[label=$\bullet$, leftmargin=*, nosep, noitemsep]
	\item Let $\Omega\subset\mathbb R^n$ be a bounded domain with boundary of class at least \(C^1\). 
	For a given time horizon $T>0$, we denote by
	\[
	Q_T=\Omega\times(0,T),\qquad \partial_lQ_T=\partial\Omega\times(0,T),
	\]
	the space-time cylinder and its lateral boundary, respectively; the parabolic boundary is
	\[
	\partial_pQ_T=\partial_lQ_T\cup(\overline\Omega\times\{0\}).
	\]
	Correspondingly, for the whole space we set
	\[
	R_T=\mathbb R^n\times(0,T),\qquad \partial_lR_T=(\mathbb{R}^n\setminus\Omega)\times(0,T),\qquad 
	\partial_pR_T=\partial_lR_T\cup(\overline\Omega\times\{0\}).
	\]
	
	\item Two families of smooth functions with compact support frequently appear:
	$C^{\infty}_0(Q_{T})$ consists of smooth functions vanishing near the lateral boundary $\partial_lQ_T$, whereas $C^{\infty}_{0,\partial p}(Q_{T})$ contains those vanishing near the full parabolic boundary $\partial_pQ_T$.
	Their closures in the $H^{1,1}$-norm yield the classical Sobolev spaces with vanishing trace:
	\[
	H^{1,1}_0(Q_{T})=\overline{C_0^{\infty}(Q_{T})}^{\|\cdot\|_{H^{1,1}(Q_{T})}},
	\qquad
	{H}^{1,1}_{0,\partial p}(Q_{T})=\overline{C^{\infty}_{0,\partial p}(Q_{T})}^{\|\cdot\|_{H^{1,1}(Q_{T})}}.
	\]
	
	
	\item For the pointwise regularity analysis we employ H\"older spaces. 
	Given $\alpha,\beta\in(0,1)$, the parabolic H\"older seminorm is
	\[
	[u]_{C^{\alpha,\beta}_{x,t}(\Omega\times (0,T))}:=
	\sup_{\begin{subarray}{c} x,x^{\prime}\in\Omega \\ t,t^{\prime}\in (0,T)\end{subarray}}
	\frac{|u(x,t)-u(x^{\prime},t^{\prime})|}{|x-x^{\prime}|^{\alpha} + |t-t^{\prime}|^{\beta}}.
	\]
	When only spatial or temporal regularity is relevant, we also use
	\[
	[u]_{C^{\alpha}_{x}(\Omega\times (0,T))}:=\sup_{t\in (0,T)}[u(\cdot,t)]_{C^{\alpha}(\Omega)},\qquad
	[u]_{C^{\beta}_{t}(\Omega\times (0,T))}:=\sup_{x\in\Omega}[u(x,\cdot)]_{C^{\beta}((0,T))}.
	\]
	The subscripts $x$ and $t$ will be omitted when no confusion can arise.
	
	\item Finally, to handle the nonlocal nature of the fractional Laplacian we introduce appropriate tail spaces. 
	The tail space is
	\[
	L^{1}_{2s}(\mathbb{R}^{n}):=\left\{v\in L^{1}_{\mathrm{loc}}(\mathbb{R}^{n}):
	\int_{\mathbb{R}^{n}}\frac{|v(x)|}{1+|x|^{n+2s}}\,dx< +\infty\right\}.
	\]
\end{itemize}

\subsection{Function spaces and weak solutions}
We define the energy space
\[
X^{1}(\mathbb R^n):=\left\{u\in L^2(\mathbb R^n): \nabla u\in L^2(\mathbb R^n),\ [u]_s^2< +\infty\right\},
\]
equipped with the norm
\[
\|u\|_{X^{1}(\mathbb R^n)}=\left(\|u\|_{L^2(\mathbb R^n)}^2+\|\nabla u\|_{L^2(\mathbb R^n)}^2+[u]_s^2\right)^{1/2},
\]
where $[u]_s$ is the Gagliardo seminorm 
\[ [u]_{s}^{2} := \iint_{\mathbb{R}^{2n}} \frac{|u(x)-u(y)|^{2}}{|x-y|^{n+2s}} dxdy .\]
We define $X^1_0(\Omega)$ is the closure of $C_{0}^{\infty}(\Omega)$ under the above norm. It is easy to see 
\[ X^1_0(\Omega) = \{ u\in H^{1}(\mathbb{R}^{n}) : u|_{\Omega}\in H_{0}^{1}(\Omega), u\equiv 0 \text{ a.e. in } \mathbb{R}^{n}\setminus\Omega \}.\]

By the classical Sobolev inequality, we can directly obtain the following lemma.

\begin{lemma}\label{Lem: space continuous embedding}
	The space $X^1_0(\Omega)$ is  continuously embedded  into $L^{p}(\Omega)$, where
	\begin{equation}\label{range of p}
		1<p< \begin{cases} 
			+\infty, & n = 1, 2, \\
			\frac{2n}{n-2}, & n > 2.
		\end{cases}
	\end{equation}
\end{lemma}

The related anisotropic parabolic energy spaces  
\[
X^{1,0}(R_{T}) := \left\{ u \in L^{2}(R_T) : 
\nabla u \in L^{2}(R_T),\ 
\int_{0}^{T} [u(\cdot,t)]_{s}^{2} \, dt < +\infty \right\},
\]
equipped with the norm
\[
\|u\|_{X^{1,0}(R_{T})} := 
\left( \|u\|_{L^{2}(R_{T})}^2 + \|\nabla u\|_{L^{2}(R_{T})}^2 + \int_{0}^{T} [u]_{s}^{2} \, dt \right)^{1/2}.
\]
Denote $ X^{1,0}_0(Q_{T})$
as the completion of $C_0^{\infty}(Q_{T})$ under this norm.  
It is easy to see $X^{1,0}_0(Q_{T})$ is isomorphic to 
\[\left\{\ \varphi\in L^{2}(0,T;H^{1}(\mathbb{R}^{n})): \varphi \text{ has spatial support compactly contained in }\Omega \  \right\}. \]

\begin{definition}\label{parabolicSobolev}
	We define the parabolic Sobolev space related to problem~\eqref{problem} by
	\[
	X^{1,1}(R_{T}) = \left\{ u \in L^{2}(R_T) : 
	D_t u \in L^{2}(R_T),\ 
	\nabla u \in L^{2}(R_T),\ 
	\text{and } \int_{0}^{T} [u(\cdot,t)]_{s}^{2}\, dt < +\infty \right\},
	\]
	equipped with the norm
	\[
	\|u\|_{X^{1,1}(R_{T})} := 
	\|u\|_{L^{2}(R_{T})} + 
	\|D_{t}u\|_{L^{2}(R_{T})} + 
	\|\nabla u\|_{L^{2}(R_{T})} + 
	\left(\int_{0}^{T} [u]_{s}^{2} dt\right)^{1/2}.
	\]
	Denote $ X^{1,1}_0(Q_T)$ as the closure of $C^\infty_0(Q_T)$ under the above norm.
\end{definition}

We remark that, by the continuous embedding $H^{1}(\mathbb{R}^{n}) \hookrightarrow H^{s}(\mathbb{R}^{n})$ (see \cite[Proposition 2.2]{DNPV12}) and the classical Poincaré inequality, we have the isomorphism \[X^{1,1}_0(Q_{T}) \cong H^{1,1}_0(Q_{T}).\]

\begin{definition}[Weak solution]\label{Def:weak sol}
	A function $u\in X^{1,1}_0(Q_{T})$ is  a weak solution to equation \eqref{problem}-\eqref{bdry-ext} if 
	\begin{equation}\label{weak sol equation}
		\int_{0}^{T} \int_{\Omega} (u_{t}\varphi + \nabla u\cdot\nabla\varphi) ~dxdt + \int_{0}^{T} \mathcal{E}(u,\varphi,t) ~dt = \int_{0}^{T} \int_{\Omega} f\varphi~ dxdt
	\end{equation}
	holds for all $\varphi\in X^{1,0}_0(Q_{T})$, where
	\[ 
	\mathcal{E}(u,\varphi,t) := 
	\iint_{\mathbb{R}^{n}\times\mathbb{R}^{n}} \frac{(u(x,t)-u(y,t))(\varphi(x,t)-\varphi(y,t))}{|x-y|^{n+2s}} dxdy. \]
	If additionally $u \in X^{1,1}_0(Q_{T})$ with  
	\begin{equation}\label{trace condition}
		\gamma u(x,0)=u_0 \quad \text{ a.e. in } \Omega,
	\end{equation}
	then $u$ is a weak solution to problem~\eqref{problem}-\eqref{bdry-initial}.
	Here  $\gamma u\in L^{2}(\partial_{p}\Omega)$ is the trace function satisfying
	\[ \lim_{m\to\infty} \int_{\partial_{p}\Omega} |u_{m}-\gamma u|^{2} d\sigma = 0 ,\]
	for any sequence $u_{m}\in C^{1}(\overline{\Omega})$ converging to $u$ in $H^{1}(\Omega)$.
\end{definition}

We note that, in view of Trace Theorem (see for example \cite{LzSUc68}), $\gamma u$ is uniquely determined, and thus the equality~\eqref{trace condition} makes sense. Moreover, 
by \cite[section 5.9.2, theorem 3]{Eva10}, a function $u$ in $X^{1,1}_0 (Q_{T})$ belongs to $ C([0,T];L^{2}(\Omega))$.


\subsection{A weak maximum principle}
For the mixed parabolic operator, we introduce the following weak maximum principle, which guarantees the existence of positive classical solutions of \eqref{eq:globalintime} in Section~\ref{sec: Asymptotic Radial Symmetry of Classical Solutions}.

%

\begin{theorem}[Weak maximum principle] \label{thm:max_principle}
	Let $u \in C^{2,1}(Q_T) \cap C(\overline{Q}_T)$ satisfy $u(\cdot, t) = 0$ on $\mathbb{R}^n \setminus \Omega$.
	Suppose that 
	\[
	\partial_t u - \Delta u + (-\Delta)^s u = f \leq 0 \quad (\text{or } f \geq 0) \quad \text{ in } Q_T.
	\]
	Then
	\[
	\sup_{Q_T} u(x,t) \leq \sup_{\partial_p Q_T} u_+(x,t) 
	\quad \bigl(\text{or } \inf_{Q_T} u(x,t) \geq \inf_{\partial_p Q_T} u_-(x,t) \bigr).
	\]
\end{theorem}
\begin{proof}
	We first prove the case $f < 0$ by contradiction.
	Suppose the conclusion fails. Then there exists a point $(x_0, t_0) \in \overline{Q}_T \setminus \partial_p Q_T$ such that
	\[
	u(x_0, t_0) = \max_{\overline{Q}_T} u(x,t) > 0, \quad t_0 \in (0,T].
	\]
	(If $\max_{\overline{Q}_T} u \leq 0$, the conclusion holds trivially.)
	
	At the maximum point $(x_0, t_0)$, we have
	\[
	\partial_t u(x_0, t_0) \geq 0, \qquad \Delta u(x_0, t_0) \leq 0.
	\]
	Moreover, by the definition of the fractional Laplacian and the fact that $u(x_0, t_0)$ is the global maximum (note that $u = 0$ on $\mathbb{R}^n \setminus \Omega$), we obtain
	\[
	(-\Delta)^s u(x_0, t_0)
	= c_{n,s} \, \mathrm{P.V.} \int_{\mathbb{R}^n} \frac{u(x_0, t_0) - u(y, t_0)}{|x_0 - y|^{n+2s}} \, dy \geq 0.
	\]
	Consequently,
	\[
	\partial_t u(x_0, t_0) - \Delta u(x_0, t_0) + (-\Delta)^s u(x_0, t_0) \geq 0,
	\]
	which contradicts the assumption $f(x_0, t_0) < 0$. Hence the conclusion holds for $f < 0$.
	
	Now consider the general case $f \leq 0$. For any $\varepsilon > 0$, define $h(t) = e^{-t}$ and set
	\[
	u_\varepsilon(x,t) = u(x,t) + \varepsilon h(t).
	\]
	Then in $Q_T$,
	\[
	(\partial_t - \Delta + (-\Delta)^s) u_\varepsilon
	= f + \varepsilon h'(t) = f - \varepsilon e^{-t} < 0.
	\]
	Applying the result already proved for $f < 0$ to $u_\varepsilon$, we obtain
	\(
	\sup_{Q_T} u_\varepsilon \leq \sup_{\partial_p Q_T} (u_\varepsilon)_+.
	\)
	Therefore,
	\[
	\sup_{Q_T} u \leq \sup_{Q_T} u_\varepsilon
	\leq \sup_{\partial_p Q_T} (u_\varepsilon)_+
	\leq \sup_{\partial_p Q_T} u_+ + \varepsilon \sup_{[0,T]} h(t).
	\]
	Letting $\varepsilon \to 0^+$ yields
	\[
	\sup_{Q_T} u(x,t) \leq \sup_{\partial_p Q_T} u_+(x,t).
	\]
	
	For the case $f \geq 0$, considering $-u$ and its corresponding equation gives, by a similar argument,
	\(
	\inf_{Q_T} u(x,t) \geq \inf_{\partial_p Q_T} u_-(x,t).
	\)
\end{proof}

We remark that  such weak maximum principle is also valid when $T= +\infty$.

\section{Existence and Regularity}\label{sec: Regularity of weak solutions}
In this section, we aim to prove Theorems~\ref{Sobolevregularity} - \ref{thm: reg up to bound} and Corollary~\ref{cor:global-exist-reg}.
Indeed, Subsections~\ref{sec:Existence and uniqueness} - \ref{Subsec: Global Boundedness Estimates for Weak Solutions} are devoted to showing Theorem~\ref{Sobolevregularity}\footnote{Actually, the statements of Theorem~\ref{Sobolevregularity} are divided in to the ones of Theorem~\ref{thm:existence} and Theorem~\ref{thm:unique} (in Subsection~\ref{sec:Existence and uniqueness}), Theorem~\ref{thm:interior-H21} (in Subsection~\ref{sec:Interior Sobolev regularity}),  and Theorem~\ref{Thm:global boundedness} (in Subsection~\ref{Subsec: Global Boundedness Estimates for Weak Solutions}).}.
We will prove Theorem~\ref{thm: reg up to bound} and Corollary~\ref{cor:global-exist-reg} in Subsection~\ref{subsec: reg-2holder} and Subsection~\ref{sec:Global-in-time existence} respectively.

\subsection{Existence and uniqueness}\label{sec:Existence and uniqueness}
We first summarise the properties of the operator
	\begin{align*}
	\mathcal{A}=\mathcal{L}^{-1}:L^{2}(\Omega) &\to X^1_0(\Omega)\subset L^{2}(\mathbb{R}^{n}) \\
	f &\mapsto \mathcal Af
\end{align*}
defined by $\mathcal{A}f = u$, where $u \in X^1_0(\Omega)$ is the unique weak solution of elliptic problem 
\begin{equation}\label{problem: elliptic}
	\left\{
	\begin{array}{ll}
		-\Delta u + (-\Delta)^{s} u = f & \text{in } \Omega, \\
		u = 0 & \text{in } \mathbb{R}^n \setminus \Omega.
	\end{array}
	\right.
\end{equation}

\begin{proposition}\label{Prop: operator A}
	$\mathcal A$ is a bounded, self-adjoint and compact operator on $L^{2}(\Omega)$.
\end{proposition}

\begin{proof}
	By the Fourier representation $(-\Delta+(-\Delta)^{s})^{-1}=(|\xi|^{2}+|\xi|^{2s})\mathcal F$ (see \cite[Proposition 3.3]{DNPV12}), $\mathcal A$ is self-adjoint.  
	The $L^{2}$-theory for the elliptic problem \cite[Theorem 1.1]{BDVa22} gives
	\[
	\|\mathcal A f\|_{X^{1}(\Omega)}\le C\|f\|_{L^{2}(\Omega)},
	\]
	so $\mathcal A$ is bounded.
	
	The compactness follows from the Rellich--Kondrachov theorem or the spectral theory of $-\Delta$ and $(-\Delta)^{s}$ (see \cite{Eva10} and \cite[Proposition 9]{SV13} respectively).
\end{proof}

Consequently, the Hilbert-Schmidt theorem yields an orthonormal basis $\{\omega_{i}\}_{i = 1}^{\infty}$ in $L^{2}(\Omega) \cap H_0^1(\Omega)$
consisting of eigenfunctions of $\mathcal A$:
\[
\mathcal A\omega_{i}=\lambda_{i}\omega_{i}\ \text{ in }\Omega,\qquad\omega_{i}=0\ \text{ in }\mathbb R^{n}\setminus\Omega,\qquad\lambda_{i}>0.
\]
Thus, $\mathcal{L}\omega_{i}=\lambda_{i}^{-1}\omega_{i}$ in $\Omega$. 

\begin{theorem}[Existence of weak solution]\label{thm:existence}
	For every $f\in L^{2}(Q_{T})$ and $u_{0}\in L^{2}(\Omega)$ problem \eqref{problem}--\eqref{bdry-initial} admits a weak solution $u\in X^{1,1}_0(Q_{T})$.
\end{theorem}

\begin{proof}
	We proceed by steps.
	
	\textbf{Step 1.} \emph{Galerkin approximation.}  
	Expand $u_{0}=\sum_{i=1}^{\infty}c_{i}\omega_{i}$ and set
	\[
	u_{m}(x,t)=\sum_{i=1}^{m}c_{i}^{m}(t)\omega_{i}(x),\qquad c_{i}^{m}(0)=c_{i}
	\]
	satisfying
	\begin{equation}\label{3}
		(\partial_{t}u_{m},\omega_{k})+B[u_{m},\omega_{k};t]=(f,\omega_{k}),\qquad 1\le k\le m,
	\end{equation}
	with
	\[
	B[u,v;t]:=\int_{\Omega}\nabla u\cdot\nabla v\,dx+\mathcal E(u,v,t).
	\]
	Note that
	\[
	(\partial_t u_m,\omega_k)=\sum_{i = 1}^{m}\frac{d}{dt}c_i^m(t)(\omega_i,\omega_k)=\frac{d}{dt}c_k^m(t),
	\]
	and
	\[
	B[u_{m},\omega_{k};t]=\sum_{i = 1}^{m}c_i^m(t)((-\Delta + (-\Delta)^{s})\omega_i,\omega_k)=\frac{1}{\lambda_k}c_k^m(t)
	\]
	via Proposition~\ref{Prop: operator A}, we obtain the linear ordinary differential equation
	\[
	\frac{d}{dt}c_k^m(t)=-\frac{1}{\lambda_k}c_k^m(t)+f_k(t)\qquad f_{k}(t)=(f(\cdot,t),\omega_{k}).
	\]
	Hence
	\[
	c_k^m(t)=e^{-t/\lambda_k}\left(c_k+\int_{0}^{t}e^{\tau/\lambda_k}f_k(\tau)d\tau\right).
	\]
	
	\textbf{Step 2.} \emph{Energy estimates.}  
	Multiply the equation~\eqref{3} by $c_{k}^{m}(t)$ and sum over $k$, we obtain
	\[
	(\partial_t u_m,u_m)+B[u_m,u_m;t]=(f,u_m),
	\]
	that is
	\[
	\frac{1}{2}\frac{d}{dt}\|u_m(\cdot,t)\|_{L^2(\Omega)}^2=-B[u_m,u_m;t]+(f(\cdot,t),u_m(\cdot,t)).
	\]
	Integrating in time, we further have
	\begin{align*}
		\frac{1}{2}\|u_m(\cdot,t)\|_{L^2(\Omega)}^2-\frac{1}{2}\|u_m(\cdot,0)\|_{L^2(\Omega)}^2
		=-\int_0^tB[u_m,u_m;s]ds+\iint_{Q_t}fu_mdxdt.
	\end{align*}
	Using the Cauchy inequality with \(\varepsilon\) and the Poincar\'e inequality, we obtain
	\begin{equation}\label{4}
		\begin{aligned}
			\sup_{t\in[0,T]}\|u_m(\cdot,t)\|_{L^2(\Omega)}^2+\int_0^TB[u_m,u_m;t]dt
			\leq C\left(\|u_m(\cdot,0)\|_{L^2(\Omega)}^2+\iint_{Q_T}f^2dxdt\right),
		\end{aligned}
	\end{equation}
	where \(C\) is a constant independent of \(m\). 
	
	Multiply the equation \eqref{3} by \(\frac{d}{dt}c_k^m(t)\) and then sum up.
	Integrating over \((0,T)\) and using the Cauchy inequality, we have
	\begin{equation}\label{5}
		\begin{aligned}
			\iint_{Q_T}|\partial_t u_m|^2dxdt+B[u_m,u_m;T]
			\leq B[u_m,u_m;0]+\iint_{Q_T}f^2dxdt.
		\end{aligned}
	\end{equation}
	
	Combining \eqref{4} and \eqref{5}, we get
	\begin{equation}\label{6}
		\iint_{Q_T}\big(|u_m|^2+|\nabla u_m|^2+|\partial_t u_m|^2\big)dxdt+\int_0^T[u_m]_s^2dt\leq C,
	\end{equation}
	where \(C >0\) is a constant independent of \(m\).
	
	\textbf{Step 3.} \emph{Passage to the limit.}  
	According to the energy estimates \eqref{6}, we see that \(\{u_m\}\) is bounded in \(X^{1,1}_0(Q_T)\). Consequently, there exists a subsequence \(\{u_{m_l}\}_{l = 1}^{+\infty}\) converges weakly to a limit \(u\in X^{1,1}_0(Q_T)\), such that
	\[
	u_{m_l}\rightharpoonup u\,\text{ in }L^2(0,T;X^1_0(\Omega)),\qquad
	\partial_t u_{m_l}\rightharpoonup u\,\text{ in }L^2(Q_T).
	\] 
	It is easy to see $u$ satisfies the integral equation \eqref{weak sol equation} for every test function of the form 
	\begin{equation}\label{7}
		v(t)=\sum_{k = 1}^{N}d^{k}(t)\omega_k,\text{ where }\{d^k(t)\}\text{ are given smooth functions},
	\end{equation}
	and hence for every $v\in L^{2}(0,T;X^1_0(\Omega))$, as functions of the form \eqref{7} are dense in this space.
	In fact, we can multiply \eqref{3} by \(d^k(t)\), sum \(k = 1,\cdots,N\), and then integrate with respect to \(t\) to find
	\[\int_0^T\big(\partial_t u_m,v\big)dt+\int_0^TB[u_m,v;t]dt=\int_0^T(f,v)dt.\]
	Set \(m = m_l\). We obtain upon passing to weak limits that $u$ satisfies the integral equation \eqref{weak sol equation} for $v$.
	
	We then verify \(\gamma u(x,0)=u_0(x)\). Note that
	\begin{align*}
		&\int_{\Omega}|u(x,t)-u_0(x)|^2dx\\
		\leq& \int_{\Omega}|u(x,t)-u_m(x,t)|^2dx+\int_{\Omega}\big|\sum_{i = 1}^{m}(c_i^m(t)-c_i)\omega_i(x)\big|^2dx+\int_{\Omega}\big|\sum_{i = m + 1}^{+\infty}c_i\omega_i(x)\big|^2dx\\
		=&:I_1 + I_2+I_3.
	\end{align*}
	When \(m\) is large enough, \(I_1\) and \(I_3\) can be arbitrarily small, and \(c_i^m(0)=c_i\). 
	
	Thus, $u$ is a weak solution to problem \eqref{problem} - \eqref{bdry-initial}.
\end{proof}

Further, we can obtain the uniqueness of weak solution by the following lemma.

\begin{lemma}[Leibniz rule with trace]\label{Lem: Leibniz Integral Rule (trace)}
	For every $u\in X^{1,1}_0(Q_T)$ and a.e. $t\in(0,T)$,
	\[
	\int_\Omega u^2(x,t)\,dx-\int_\Omega(\gamma u(x,0))^2\,dx=2\int_0^t\int_\Omega u\partial_t u\,dxds.
	\]
\end{lemma}

\begin{proof}[Proof sketch]
	Approximate $u$ by a sequence in $ C^\infty_0(Q_T)$ and pass to the limit.
\end{proof}

\begin{theorem}[Uniqueness of weak solution]\label{thm:unique}
	Under the assumptions of Theorem \ref{thm:existence}, the weak solution of \eqref{problem}--\eqref{bdry-initial} is unique.
\end{theorem}

\begin{proof}
	Let $u=u_1-u_2$ be the difference of two weak solutions. Then $u\in X^{1,1}_0(Q_T)$, $\gamma u(\cdot,0)=0$.
    Choosing the test function as $\varphi=u\chi_{[0,r]}$, $r\in[0,T]$, we have
	\[
	\iint_{Q_r}(u_tu + |\nabla u|^2)dxdt+\int_0^r[u]_s^2dt = 0.
	\]
	Thus
	$
	\iint_{Q_r}uu_tdxdt\leqslant 0.
	$
	Using Lemma \ref{Lem: Leibniz Integral Rule (trace)}, we deduce
	\[
	\int_{\Omega}u^2(x,r)dx\leqslant 0, \quad \text{ a.e. } r\in[0,T].
	\]
	Hence, \(u = 0\) a.e. in \(Q_T\).
\end{proof}

\subsection{Interior Sobolev regularity of weak solution}\label{sec:Interior Sobolev regularity}
Having established the existence and uniqueness of weak solutions, we now turn to their interior smoothness.
%
%
%

\begin{theorem}[Interior $H^{2,1}$-regularity]\label{thm:interior-H21}
	Let $f\in L^2(Q_T)$ and $u\in X^{1,1}_0(Q_T)$ be the weak solution of \eqref{problem}--\eqref{bdry-initial}.  
	For $\delta>0$, set $\Omega_{\delta}:= \{ x\in \Omega; \text{dist}(x, \partial \Omega) > \delta\}$ and $ Q_{T}^{\delta}= \Omega_{\delta} \times (0,T).$
	Then $u\in H^{2,1}(Q_T^\delta)$ and
	\begin{equation}\label{Interior Sobolev estimates}
	\|u\|_{H^{2,1}(Q_T^\delta)}\le C\bigl(\|u\|_{H^{1,0}(Q_T)}+\|f\|_{L^2(Q_T)}+\|u_0\|_{H^1(\Omega)}\bigr)
	\end{equation}
	with some constant $C=C(n, \delta) > 0$.
\end{theorem}

\begin{proof}
	We proceed by steps.
	
	\textbf{Step 1.} \emph{Difference-quotient argument with truncation method.}  
	Choose  
	\[ \zeta(x) \in C_{0}^{\infty}(\mathbb{R}^{n})\ \text{ with }  \text{supp}(\zeta)\subset \Omega, \ 0 \leq \zeta \leq 1, \ \zeta \equiv 1 \text{ on } \Omega_{\delta},\ |\nabla \zeta| \leq C/\delta.\]
	
	For $|h|<\frac12\operatorname{dist}(\operatorname{supp}\zeta,\partial\Omega)$ insert the test function  
	\[
	\varphi=D_k^{-h}(\zeta^2 D_k^h u)\in X^{1,0}_0(Q_T)
	\]  
	into the weak formulation \eqref{weak sol equation}.  
	The identity  
	\begin{align*}
		&\iint_{Q_T} u_t\varphi\,dxdt
		=-\iint_{Q_T} D_k^{-h}u_t \cdot \zeta^{2} D_{k}^{h} u\ dxdt\\
		=&\iint_{Q_T} \frac{\partial}{\partial t}D_k^hu \cdot \zeta^{2} D_{k}^{h} u\ dxdt=\frac{1}{2}\iint_{Q_T}\frac{\partial}{\partial t}(\zeta^{2} (D_{k}^{h}u)^2) \ dxdt\\
		=&\frac{1}{2}\int_{\Omega}\zeta^{2} |D_{k}^{h}u(x,T)|^2 dx-\frac{1}{2}\int_{\Omega}\zeta^{2} |D_{k}^{h}u(x,0)|^2 dx
	\end{align*}
	yields 
	\begin{align*}
		&\iint_{Q_T}  \nabla u\cdot\nabla(D_{k}^{-h}(\zeta^{2} D_{k}^{h} u)) dxdt + \int_{0}^{T} \mathcal{E}(u,D_{k}^{-h}(\zeta^{2} D_{k}^{h} u),t) dt\\ 
		\leq&\iint_{Q_T} fD_{k}^{-h}(\zeta^{2} D_{k}^{h} u) dxdt+\frac{1}{2}\int_{\Omega}\zeta^{2} (D_{k}^{h}u_0(x))^2 dx.
	\end{align*}
	Combining this with the elliptic estimate \cite[(4.3)]{BDVa22}, we have 
	\[
	\iint_{Q_T^\delta}\big|D_k^h \ \nabla u\big|^2\,dxdt\le C(n,\delta)\left( \|u\|_{L^2(0,T;H^1(\Omega))}^2+\|f\|_{L^2(Q_T)}^2+ \|D_k^hu_0(x)\|_{L^2(\Omega)}^2\right).
	\]  
	Using \cite[Theorem 3 in section 5.8.2]{Eva10}, and letting \(h \to 0\), we have  
	\begin{equation}\label{11}
		\iint_{Q_T^\delta}|D^2 u|^2\,dxdt\le C(n,\delta)\left(\|u\|_{L^2(0,T;H^1(\Omega))}^2+\|f\|_{L^2(Q_T)}^2+\|u_0\|_{H^1(\Omega)}^2\right).
	\end{equation}
	
	\textbf{Step 2.} Time-derivative bound. 
	From the equation $u_t=\Delta u-(-\Delta)^s u+f$, we obtain
	\begin{equation*}\label{9}
		\begin{aligned}
			\iint_{Q_T^{\delta}}| \partial_t u|^2dxdt 
			\leq C(n,\delta) \left( \iint_{Q_T^{\delta}}|\Delta u|^2dxdt+\iint_{Q_T^{\delta}}\left| (-\Delta)^su\right|^2dxdt+\iint_{Q_T^{\delta}}|f|^2dxdt\right).
		\end{aligned}
	\end{equation*}
	
	Since $(-\Delta)^s$ maps the Bessel potential space $H^{p}(\mathbb R^n)$ continuously into $H^{p-2s}(\mathbb R^n)$ for any $p\in\mathbb R$ (see, e.g., \cite{Ste70}), we have 
	\begin{equation}\label{66}
		\|(-\Delta)^s u\|_{L^2(\mathbb R^n)}\le C(n,s) \|u\|_{H^2(\mathbb R^n)}.
	\end{equation}
	
	Combining this with \eqref{11} and \eqref{66} we obtain  
	\[
	\|\partial_t u\|_{L^2(Q_T^\delta)}\le C(n,\delta)\bigl(\|D^2 u\|_{L^2(Q_T^\delta)}+\|u\|_{L^2(0,T;H^2(\Omega))}+\|f\|_{L^2(Q_T)}\bigr).
	\]
	
	Hence, by combining the estimates in both steps, we  obtain \eqref{Interior Sobolev estimates}.
\end{proof}

One may have the following higher-order interior regularity result.
\begin{theorem}[Higher-order interior regularity]\label{thm:high-int}
	Assume $\partial\Omega\in C^\infty$, $u_0\in H^{2k+1}(\Omega)$, $f\in H^{2k,k}(Q_T)$ and the compatibility conditions  
	\[
	g_0:=u_0\in H^1_0(\Omega),\quad g_j:=\partial_t^{j-1}f(\cdot,0)-\mathcal L g_{j-1}\in H^1_0(\Omega),\quad 1\le j\le k.
	\]  
	Then $u\in H^{2k+2,k+1}(Q_T^\delta)$ and  
	\[
	\|u\|_{H^{2k+2,k+1}(Q_T^\delta)}\le C\bigl(\|u\|_{H^{1,0}(Q_T)}+\|f\|_{H^{2k,k}(Q_T)}+\|u_0\|_{H^{2k+1}(\Omega)}\bigr)
	\]
	with some positive constant \(C=C(n,\delta)\).
\end{theorem}

\begin{proof}[Sketch of proof]
	Following \cite[Theorem 6, $\sec$ 7.1]{Eva10}, we freeze time and view the parabolic equation as an elliptic problem for $u(\cdot,t)$ with right-hand side $f(\cdot,t)-\partial_t u(\cdot,t)$.  
	Applying the interior elliptic regularity \cite[Theorem 1.4]{BDVa22} gives
	\[
	\|u(\cdot,t)\|_{H^{2k+2}(\Omega_\delta)}\le C\bigl(\|f(\cdot,t)\|_{H^{2k}(\Omega)}+\|\partial_t u(\cdot,t)\|_{H^{2k}(\Omega)}\bigr),
	\]
	where \(\Omega_{\delta}:= \{ x\in \Omega; \text{dist}(x, \partial \Omega) > \delta\} \).
	An induction on $k$ together with the energy estimates for $\partial_t^k u$ yields the claimed space-time regularity.
\end{proof}

\subsection{Global Boundedness of weak solution}\label{Subsec: Global Boundedness Estimates for Weak Solutions} 
In this subsection, we prove that the weak solution is globally bounded.  
To begin with, we introduce the following notion of essential supremum/infimum for functions in $ X^{1,1}_0(Q_T)$.

\begin{definition}[Supremum and infimum]
	For $u\in X^{1,1}_0(Q_T)$ set
	\[
	\sup_{Q_T}u=\inf\bigl\{ ~l\in\mathbb R:(u-l)_+=0\ \text{a.e. in }Q_T~ \bigr\}, \qquad
	\sup_{\partial_pQ_T}u=\inf\bigl\{~l\in\mathbb R:(u-l)_+\in H^{1,1}_{0,\partial p}(Q_T)~\bigr\},
	\]
	and define $\inf_{Q_T}u=-\sup_{Q_T}(-u)$, $\inf_{\partial_pQ_T}u=-\sup_{\partial_pQ_T}(-u)$.
\end{definition}

\begin{lemma}\label{Lem: v,v_+}
	For every $v\in X^{1,1}(\mathbb R^n\times(0,T))$,
	\[
	(v(x,t)-v(y,t))(v_+(x,t)-v_+(y,t))\ge|v_+(x,t)-v_+(y,t)|^2\quad\forall x,y\in\mathbb R^n.
	\]
\end{lemma}
\begin{proof}[Proof sketch]
	Substitute $v = v_+ - v_-$ and expand.
\end{proof}

We now give a well known basic iteration lemma. (See \cite{WS15} for a proof.)
\begin{lemma}[Basic Iteration Lemma]\label{Lem: Basic Iteration Lemma}
	Let $\varphi(t)$ be a nonnegative and nonincreasing function on $[k_0, +\infty)$ satisfying
	\[
	\varphi(h) \leq \left( \frac{M}{h - k} \right)^{\alpha} [\varphi(k)]^{\beta}, \quad \forall h > k \geq k_0
	\]
	for constants $M > 0$, $\alpha > 0$, $\beta > 1$. Then
	\[
	\varphi(k_0+d)=0,
	\]
	where $d=M2^{\beta/(\beta-1)}[\varphi(k_0)]^{(\beta-1)/\alpha}$.
\end{lemma}

\begin{theorem}[Weak maximum principle for weak solution]\label{Thm: Weak maximum principle}
	Let $f\in L^\infty(Q_T)$ and $u\in X^{1,1}_0(Q_T)$ be a weak solution of \eqref{problem}. Then there exists a constant $C = C(n, \Omega,T)$ such that
	\[
	\sup_{Q_T}u\le\sup_{\partial_pQ_T}u+C\|f\|_{L^\infty(Q_T)}.
	\]
\end{theorem}
\begin{proof}
	Let $l := \sup_{\partial_{p}Q_{T}} u$. Fix an arbitrary level \(k>l\). Then \(u-k\) is still a weak solution.
	Choose the test function  
	\[
	\varphi = (u-k)_{+} \chi_{[t_{1}, t_{2}]},
	\]  
	where \(\chi_{[t_{1},t_{2}]}\) denotes the characteristic function of the interval \([t_{1},t_{2}]\subset[0,T]\).  Substituting \(\varphi\) into the weak formulation yields
	\begin{equation}
		\begin{aligned}\label{uk weak formula}
			&\underbrace{\int_{0}^{T} \int_{\Omega} (u-k)_t (u-k)_+ \chi_{[t_{1}, t_{2}]} \, dx \, dt}_{\text{I}} 
			+ \underbrace{\int_{0}^{T} \int_{\Omega} \chi_{[t_{1}, t_{2}]} |\nabla (u-k)_+|^2 \, dx \, dt}_{\text{II}} \\
			&+ \underbrace{\int_{0}^{T} \int_{\mathbb{R}^{2n}} \bigl( (u(x,t)-k) - (u(y,t)-k) \bigr) \bigl( (u(x,t)-k)_+ - (u(y,t)-k)_+ \bigr) \, d\mu \, dt}_{\text{III}} \\
			&= \int_{0}^{T} \int_{\Omega} f (u-k)_{+} \chi_{[t_{1}, t_{2}]} \, dx \, dt.
		\end{aligned}
	\end{equation}
	Here and in the sequel, $d\mu = |x-y|^{-n-2s}dxdy$ denotes the standard fractional measure.
	
	By Lemma \ref{Lem: v,v_+}, the nonlocal term satisfies $\text{III} \geq 0$. Moreover, 
	\[
	\text{I} = \frac{1}{2} \int_{t_1}^{t_2} \frac{d}{dt} \int_{\Omega} |(u-k)_{+}|^{2} \, dx \, dt \geq 0.
	\]
	Dropping the non-negative summands \(\mathrm{I}\) and \(\mathrm{III}\) in \eqref{uk weak formula} we arrive at
	\[
	\text{II} \leq \int_{t_1}^{t_2} \int_{\Omega} f (u-k)_{+} \, dx \, dt.
	\]
	The rest of the proof follows the standard De-Giorgi iteration for local \(L^{\infty}\)-bounds.  For the reader's convenience we include the details.
	
	Let $I_{k}(t) = \int_{\Omega} (u-k)_{+}^{2} \, dx$, and let $\sigma$ be a maximum point of $I_{k}(t)$ on $[0,T]$. Since $I_{k}(0) = 0$, we may assume $\sigma > 0$. For any sufficiently small $\varepsilon > 0$, take $t_1 = \sigma - \varepsilon$ and $t_2 = \sigma$. Then
	\[
	\int_{\sigma - \varepsilon}^{\sigma} \int_{\Omega} |\nabla (u-k)_+|^2 \, dx \, dt \leq \int_{\sigma - \varepsilon}^{\sigma} \int_{\Omega} |f| (u-k)_{+} \, dx \, dt.
	\]
	Dividing both sides by $\varepsilon$ and letting $\varepsilon \to 0^{+}$, we obtain
	\begin{equation}\label{1}
		\int_{\Omega} |\nabla (u(x, \sigma)-k)_+|^2 \, dx \leq \int_{\Omega} |f(x,\sigma)| (u(x, \sigma)-k)_{+} \, dx.
	\end{equation}
	
	Define
	$
	A_{k}(t) = \{ x \mid u(x,t) > k \}, \, \mu_{k} = \sup_{0 < t < T} |A_{k}(t)|.
	$ 
	Since $X^1_0(\Omega) \hookrightarrow L^{p}(\Omega)$, where $p$ is defined as in \eqref{range of p}, we have
	\[
	\left( \int_{A_k(\sigma)} (u - k)_+^p \, dx \right)^{1/p} \leq \left( \int_{\Omega} |\nabla(u(x, \sigma) - k)_+|^2 \, dx \right)^{1/2}.
	\]
	
	Combining this with \eqref{1} and setting $q = \frac{p}{p-1}$ yields
	\begin{align*}
		\left( \int_{A_k(\sigma)} (u - k)_+^p \, dx \right)^{2/p}
		\leq C \left( \int_{A_k} |f(x,\sigma)|^{q} \, dx \right)^{1/q} \left( \int_{A_k} (u(x, \sigma)-k)_{+}^{p} \, dx \right)^{1/p}.
	\end{align*}
	Consequently,
	\begin{align*}
		\left( \int_{A_k(\sigma)} (u - k)_+^p \, dx \right)^{\frac{1}{p}}
		\leq C \left( \int_{A_k} |f(x,\sigma)|^{q} \, dx \right)^{\frac{1}{q}} \leq C \|f\|_{L^\infty(Q_T)} |A_k(\sigma)|^{\frac{1}{q}} \leq C \|f\|_{L^\infty(Q_T)} \mu_k^{\frac{1}{q}}.
	\end{align*}
	
	Applying H\"older inequality to $I_k(\sigma)$, we obtain
	\begin{align*}
		I_k(\sigma) \leq \left( \int_{A_k(\sigma)} (u - k)_+^p \, dx \right)^{2/p} |A_k(\sigma)|^{(p-2)/p} \leq (C \|f\|_{L^\infty(Q_T)})^{2} \mu_k^{(3p-4)/p}.
	\end{align*}
	Therefore, for any $h > k$ and $t \in [0, T]$,
	\[
	(h - k)^2 \mu_h \leq I_k(t) \leq I_k(\sigma) \leq (C \|f\|_{L^\infty(Q_T)})^{2} \mu_{k}^{(3p-4)/p}.
	\]
	That is,
	\[
	\mu_h \leq \left( \frac{C\|f\|_{L^\infty(Q_T)}}{h - k} \right)^2 \mu_k^{(3p-4)/p}, \quad \text{ where } \frac{3p - 4}{p} > 1.
	\]
	
	Applying Lemma \ref{Lem: Basic Iteration Lemma}, we conclude that
	\[
	\mu_{l+d} = \sup_{0 < t < T} |A_{l+d}(t)| = 0,
	\]
	where
	$
	d = C\|f\|_{L^\infty(Q_T)} \mu_l^{1-2/p} 2^{(3p-4)/(2p-4)} \leq C|\Omega|^{1-2/p} 2^{(3p-4)/(2p-4)} \|f\|_{L^\infty(Q_T)}.
	$
	
	By the definition of $A(k)$, we conclude that 
	\[
	\sup_{Q_T}u\le\sup_{\partial_pQ_T}u+C\|f\|_{L^\infty(Q_T)}.\qedhere
	\]
\end{proof}

As a direct corollary of Theorem \ref{Thm: Weak maximum principle}, we obtain the following result:
\begin{theorem}[Global boundedness]\label{Thm:global boundedness}
	Let $f \in L^{\infty}(Q_{T})$ and suppose $u \in X^{1,1}_0(Q_{T})$ is a weak solution to \eqref{problem}-\eqref{bdry-initial}. Then $u \in L^{\infty}(Q_{T})$ and 
	\begin{equation*}
		\|u\|_{L^\infty ( Q_{T})} \leq \|u_0\|_{L^\infty ( \Omega)} + C \|f\|_{L^{\infty}( Q_{T})}.
	\end{equation*}
	for some constant $C = C(n, \Omega,T)>0$.
\end{theorem}

\subsection{$C^{2+\alpha,1+\alpha/2}$-regularity up to the boundary}\label{subsec: reg-2holder}
We establish the H\"older regularity  
\(
u\in C^{2+\alpha,1+\alpha/2}(\overline Q_T)
\) 
for the mixed parabolic operator \(\mathcal L^t=\partial_t-\Delta+(-\Delta)^s\) by means of the continuity method.  The key observation is that the nonlocal term $(-\Delta)^s$ can be treated as a lower-order perturbation of the classical heat operator $\partial_t-\Delta$ provided \(0<s<1/2\)\footnote{The restriction $s<1/2$ is somewhat unavoidable in view of the stationary case \cite[Remark 4.6]{CSZ25}}, and $0<\alpha\le 1-2s$.

Throughout this subsection, let $\Omega\subset\mathbb R^n$ be a bounded $C^{2+\alpha}$-domain,  
\[
0<s<1/2,\qquad 0<\alpha\le 1-2s.
\]  
and we work in the space  
\begin{equation}\label{space: B}
	\mathcal{B}_{\alpha}(Q_{T}):= \left\{ ~u\in C(R_T): u|_{Q_T} \in C^{2+\alpha,1+\frac{\alpha}{2}}(\overline Q_T), \, u(\cdot,t)|_{ \mathbb{R}^n \setminus \Omega} \equiv 0  ,u(\cdot,0)|_{\Omega}=0\  \right\}.
\end{equation}

\subsubsection{H\"older estimates for the fractional term}
We first quantify the mapping property of $(-\Delta)^s$ on $\mathcal B_\alpha(Q_T)$.
\begin{proposition}\label{prop:frac-Holder}
	For every $u\in\mathcal B_\alpha(Q_T)$ we have $(-\Delta)^s u\in C^{\alpha,\alpha/2}_{x,t}(\overline Q_T)$ and  
	\begin{equation}\label{est:frac-Holder}
		\|(-\Delta)^s u\|_{C^{\alpha,\alpha/2}_{x,t}(\overline Q_T)}\le C\|u\|_{C^{1,1/2}_{x,t}(\overline Q_T)},
	\end{equation}
	where $C=C(n,s,\alpha,\Omega)$.
\end{proposition}
\begin{proof}
	Decompose the full H\"older norm into seminorm and supremum:  
	\[
	\|(-\Delta)^{s}u\|_{C^{\alpha,\alpha/2}_{x,t}(\overline{Q}_T)}
	= [(-\Delta)^{s}u]_{C^{\alpha,\alpha/2}_{x,t}(\overline{Q}_T)}
	+ \|(-\Delta)^{s}u\|_{L^{\infty}(\overline{Q}_T)}.
	\] 
	We bound each component separately.
	
	\textit{Supremum norm.} By \cite[Lemma 5.5]{SVWZ25}, there exists some conatnat \(0<C=C(n,s,\alpha,\Omega)\) such that
	\[\|(-\Delta)^{s}u\|_{L^{\infty}(\overline{Q}_{T})}=\sup_{t\in(0,T)}\|(-\Delta)^{s}u(\cdot,t)\|_{L^{\infty}(\overline{\Omega})}
	\le \sup_{t\in(0,T)}C\,\|u(\cdot,t)\|_{C^{1}(\overline{\Omega})}
	\le C\,\|u\|_{C^{1,1/2}_{x,t}(\overline{Q}_T)}.\]
	
	\textit{H\"older seminorm.} Rather than bounding the full parabolic H\"older seminorm directly, we invoke the spatial-temporal equivalence  
	\[
	[(-\Delta)^{s}u]_{C^{\alpha,\alpha/2}_{x,t}(\overline{Q}_T)}
	\sim
	[(-\Delta)^{s}u]_{C^{\alpha}_{x}(\overline{Q}_T)}
	+
	[(-\Delta)^{s}u]_{C^{\alpha/2}_{t}(\overline{Q}_T)}
	\]  
	and estimate each term separately.
	
	The spatial part is controlled by  
	\begin{equation}\label{63}
		[(-\Delta)^s u]_{C^{\alpha}_{x}(\overline Q_T)}
		=\sup_{t\in[0,T]}[(-\Delta)^s u(\cdot,t)]_{C^{\alpha}(\overline\Omega)}
		\le C(n,s,\alpha,\Omega)\,[u]_{C^{1}_{x}(\overline Q_T)},
	\end{equation}
	where \([u]_{C^{1}_{x}}\) denotes the Lipschitz seminorm of \(u\) in the spatial variable, thanks to \cite[Lemma~5.4]{SVWZ25}.
	
	As for the temporal regularity, we fix \(x\in\overline\Omega\) and choose arbitrary \(t_1,t_2\in(0,T)\). Set \(r:=|t_1-t_2|^{1/2}\). Then
	\begin{equation}\label{60}
		|(-\Delta)^s u(x,t_1)-(-\Delta)^s u(x,t_2)|
		\le I_1+I_2,
	\end{equation}
	where
	\begin{align*}
		I_1&=\Bigl|\int_{B_r(x)}\frac{u(x,t_1)-u(y,t_1)-(u(x,t_2)-u(y,t_2))}{|x-y|^{n+2s}}\,dy\Bigr|,\\
		I_2&=\Bigl|\int_{\mathbb{R}^n\setminus B_r(x)}\frac{u(x,t_1)-u(y,t_1)-(u(x,t_2)-u(y,t_2))}{|x-y|^{n+2s}}\,dy\Bigr|.
	\end{align*}
	
	For \(I_1\), using \(|u(x,t_i)-u(y,t_i)|\le [u]_{C^{1}_{x}(\overline Q_T)}|x-y|\) for \(i=1,2\), we obtain
	\[
	I_1\le 2[u]_{C^{1}_{x}(\overline Q_T)}\int_{B_r(x)}|x-y|^{1-n-2s}\,dy\le C(n,s)[u]_{C^{1}_{x}(\overline Q_T)}r^{1-2s}.
	\]	
	For \(I_2\), we invoke \(|u(x,t_1)-u(x,t_2)|\le [u]_{C^{1/2}_{t}(\overline Q_T)}|t_1-t_2|^{1/2}=[u]_{C^{1/2}_{t}(\overline Q_T)}r\),
	\[
	I_2\le 2[u]_{C^{1/2}_{t}(\overline Q_T)}r\int_{\mathbb{R}^n\setminus B_r(x)}|x-y|^{-n-2s}\,dy
	=C(n,s)[u]_{C^{1/2}_{t}(\overline Q_T)}r^{1-2s}.
	\]
	
	Combining the estimates for \(I_1\) and \(I_2\) with \eqref{60} and taking supremum over \(x\in\overline\Omega\) gives  
	\[
	\|(-\Delta)^s u(\cdot,t_1)-(-\Delta)^s u(\cdot,t_2)\|_{L^\infty(\overline\Omega)}
	\le C(n,s)[u]_{C^{1,1/2}_{x,t}(\overline Q_T)}|t_1-t_2|^{\frac{1-2s}{2}}.
	\]  
	Hence, for every \(\alpha\le 1-2s\),  
	\begin{equation}\label{temp-seminorm}
		[(-\Delta)^s u]_{C^{\alpha/2}_{t}(\overline Q_T)}
		\le C(n,s)[u]_{C^{1,1/2}_{x,t}(\overline Q_T)}.
	\end{equation}
	
	Combining \eqref{63} and \eqref{temp-seminorm} we arrive at the desired estimate  
	\[
	[(-\Delta)^s u]_{C^{\alpha,\alpha/2}_{x,t}(\overline Q_T)}
	\le C(n,s,\alpha,\Omega)[u]_{C^{1,1/2}_{x,t}(\overline Q_T)}.\qedhere
	\]
\end{proof}

\subsubsection{Uniform a-priori estimate for the continuity path}
For $\tau\in[0,1]$ set  
\begin{equation}\label{operator: L}
	\mathscr{L}_\tau = (1 - \tau )(\partial_t-\Delta) + \tau (\partial_t-\Delta+(-\Delta)^s) = \partial_t-\Delta + \tau (-\Delta)^s.
\end{equation}
Proposition~\ref{prop:frac-Holder} together with the classical parabolic Schauder theory implies the following uniform estimate.

\begin{proposition}\label{prop:uniform-est}
	There exists a constant $C=C(n,s,\alpha,\Omega,T)>0$ such that for every $\tau\in[0,1]$ and every $u\in\mathcal B_\alpha(Q_T)$,
	\begin{equation}\label{est:uniform}
		\|u\|_{C^{2+\alpha,1+\alpha/2}_{x,t}(\overline Q_T)}\le C\Bigl(\|\mathscr L_\tau u\|_{C^{\alpha,\alpha/2}_{x,t}(\overline Q_T)}+\sup_{\overline Q_T}|u|\Bigr).
	\end{equation}
\end{proposition}

\begin{proof}
	We first verify that \(\mathscr{L}_\tau u\in C^{\alpha,\alpha/2}_{x,t}(\overline Q_T)\).  
	The assertion follows at once from Proposition \ref{prop:frac-Holder} together with the observation that \(\partial_t u - \Delta u \in C^{\alpha,\alpha/2}_{x,t}(\overline{Q}_T)\) whenever \(u \in \mathcal{B}_\alpha(Q_T)\).
	
	Since \(\partial\Omega\) is of class \(C^{2+\alpha}\), the classical parabolic Schauder theory (cf. \cite[Theorem 6, Chapter 3]{LG96}) yields, for any solution \(u \in \mathcal{B}_\alpha(Q_T)\), 
	\begin{equation}\label{61}
		\begin{aligned}
			\|u\|_{C^{2+\alpha,1+\alpha/2}_{x,t}(\overline{Q}_T)}&\le C(n,\alpha,\Omega,T)\bigl(\|\mathscr L_\tau u-\tau(-\Delta)^s u\|_{C^{\alpha,\alpha/2}_{x,t}(\overline{Q}_T)}+\sup_{\overline{Q}_T}|u|\bigr)\\
			&\leq C(n,s,\alpha,\Omega,T)\Bigl(\|\mathscr{L}_\tau u\|_{C^{\alpha,\alpha/2}_{x,t}(\overline{Q}_T)}+\|u\|_{C^{1,1/2}_{x,t}(\overline{Q}_T)}+\sup_{\overline{Q}_T}|u|\Bigr),
		\end{aligned}
	\end{equation} 
	where the last inequality follows from Proposition~\ref{prop:frac-Holder}.
	
	Exploiting once more the smoothness of \(\partial\Omega\) and invoking the anisotropic Hölder interpolation inequality (see, e.g., \cite[Theorem 8.8.1]{Kry96}), we obtain, for every \(\varepsilon>0\), there exists a constant \(C(\varepsilon)>0\) such that
	\[
	\|u\|_{C^{1,1/2}_{x,t}(\overline{Q}_T)}\leq \varepsilon\|u\|_{C^{2+\alpha,1+\alpha/2}_{x,t}(\overline{Q}_T)}+C(\varepsilon)\sup_{\overline{Q}_T}|u|.
	\]
	Combining this with \eqref{61} and choosing \(\varepsilon\) sufficiently small, we complete the proof.
\end{proof}

\subsubsection{Method of continuity}
We can now set up the framework for the method of continuity.
\begin{lemma}
	Endow \(\mathcal{B}_{\alpha}(Q_{T})\) defined in \eqref{space: B} with the norm 
	\[
	\|u\|_{\mathcal B_{\alpha}(Q_{T})}:=\|u\|_{C^{2+\alpha,1+\alpha/2}_{x,t}(\overline Q_{T})}.
	\]
	Let \(\mathscr{L}_{\tau}\) be the linear operator defined in \eqref{operator: L}.  Then 
	
	\begin{enumerate}
		\item \(\mathcal{B}_{\alpha}(Q_{T})\) is a Banach space;
		
		\item \(\mathscr{L}_{0}\) and \(\mathscr{L}_{1}\) are bounded linear operators from \(\mathcal{B}_{\alpha}(Q_{T})\) into \(C^{\alpha,\alpha/2}_{x,t}(\overline{Q}_{T})\);
		
		\item there exists a constant \(C(n,s,\alpha,\Omega,T)>0\), independent of \(\tau\in[0,1]\), such that
		\begin{equation}\label{65}
			\|u\|_{\mathcal{B}_{\alpha}(Q_{T})}\le C\|\mathscr{L}_{\tau}u\|_{C^{\alpha,\alpha/2}_{x,t}(\overline{Q}_{T})},\quad\text{for all}\, u\in\mathcal{B}_{\alpha}(Q_{T}),\,\tau\in[0,1].
		\end{equation}
	\end{enumerate}
\end{lemma}

\begin{proof}
	The assertion (1) is immediate because $C^{2+\alpha,1+\alpha/2}(\overline Q_{T})$ is complete
	and the vanishing conditions defining $\mathcal B_{\alpha}(Q_{T})$ are closed.
	
	(2) For every \(u\in\mathcal{B}_{\alpha}(Q_{T})\), we have
	\[
	\|\mathscr{L}_{0}u\|_{C^{\alpha,\alpha/2}_{x,t}(\overline{Q}_{T})}
	=\|\partial_{t}u-\Delta u\|_{C^{\alpha,\alpha/2}_{x,t}(\overline{Q}_{T})}
	\leq C(n,\alpha,\Omega,T)\|u\|_{\mathcal{B}_{\alpha}(Q_{T})}.
	\]
	
	Combining this with Proposition \ref{prop:frac-Holder}, we obtain
	\[
	\|\mathscr{L}_{1}u\|_{C^{\alpha,\alpha/2}_{x,t}(\overline{Q}_{T})}
	\le\|\partial_{t}u-\Delta u\|_{C^{\alpha,\alpha/2}_{x,t}(\overline{Q}_{T})}
	+\|(-\Delta)^{s}u\|_{C^{\alpha,\alpha/2}_{x,t}(\overline{Q}_{T})}
	\le C(n,s,\alpha,\Omega,T)\|u\|_{\mathcal{B}_{\alpha}(Q_{T})}.
	\]
	Hence \(\mathscr{L}_{0}\) and \(\mathscr{L}_{1}\) are bounded operators from \(\mathcal{B}_{\alpha}(Q_{T})\) into \(C^{\alpha,\alpha/2}_{x,t}(\overline{Q}_{T})\).
	
	$(3)$ By the weak maximum principle (Theorem \ref{Thm: Weak maximum principle}), we derive
	\[\sup_{\overline Q_T}|u| \leq C(n,\Omega,T)\|\mathscr{L}_\tau u\|_{L^{\infty}(Q_T)}\leq C(n,\Omega,T)\|\mathscr{L}_\tau u\|_{C^{\alpha,\alpha/2}(\overline Q_T)}, \quad \forall \tau \in[0,1]. \]
	Inserting this into the a-priori estimate \eqref{est:uniform} we obtain
	\eqref{65}.
\end{proof}

\subsubsection{Proof of Theorem~\ref{thm: reg up to bound}}
\begin{proof}[Proof of Theorem~\ref{thm: reg up to bound}]
	We proceed by splitting the problem into zero-initial and non-zero-initial parts.
	
	\medskip
	\noindent
	\textit{Case 1: Zero initial datum.}  
	Assume $u_0=0$.  
	Since the operator \(\mathscr{L}_0 = \partial_t-\Delta\) is surjective from \(\mathcal{B}_\alpha(Q_T)\) to \(C^{\alpha,\alpha/2}(\overline Q_T)\), we may apply the method of continuity to the family  
	\[
	\mathscr L_\tau=\partial_t-\Delta+\tau(-\Delta)^s,\qquad \tau\in[0,1]
	\]  
	to conclude that $\mathscr L_1$ is also surjective.  
	Hence, for every $f\in C^{\alpha,\alpha/2}(\overline Q_T)$ there exists a unique solution $u\in\mathcal B_\alpha(Q_T)$ satisfying  
	\[
	\partial_t u-\Delta u+(-\Delta)^s u=f\quad \text{a.e. in }\Omega\times(0,T).
	\]  
	Moreover, by the parabolic Schauder estimate \eqref{est:uniform}, we have
	\begin{equation}\label{est: case1}
		\|u\|_{C^{2+\alpha,1+\alpha/2}(\overline Q_T)}
		\le C(n,s,\alpha,\Omega,T)\bigl(\|f\|_{C^{\alpha,\alpha/2}(\overline Q_T)}+\|u\|_{L^\infty(Q_T)}\bigr).
	\end{equation}
	
	\medskip
	\noindent
	\textit{Case 2: Non-zero initial datum.}  
	For general $u_0\in C^{2+\alpha}(\overline\Omega)$, set 
	\[\tilde{u}(x,t)=\left\{\begin{alignedat}{2}
		u(x,t)&-u_0(x) \qquad && \text{ in } \Omega \times[0,T),\\
		&0 \qquad && \text{ in } \mathbb{R}^n \setminus \Omega \times [0,T).
	\end{alignedat}
	\right.\]  
	Then $\tilde u$ solves  
	\begin{equation*}
		\left\{
		\begin{alignedat}{2}
			\partial_t \tilde{u} - \Delta \tilde{u} + (-\Delta)^s \tilde{u} &= f- \mathcal{L}u_0 \qquad &&\text{in }  \Omega \times (0,T),\\
			\tilde{u}(x,t) &= 0 \qquad && \text{in } \bigl(\mathbb R^n\setminus\Omega\bigr) \times [0,T), \\
			\tilde{u}(x,0) &= 0 \quad && \text{in } \Omega.
		\end{alignedat}
		\right.
	\end{equation*}
	Since $(-\Delta)^s u_0\in C^{\alpha}(\overline\Omega)$ by Proposition \ref{prop:frac-Holder}, it follows that
	\[
	f-\mathcal{L} u_0\in C^{\alpha,\alpha/2}(\overline Q_T).
	\]  
	Applying the result from \textit{Case 1} to $\tilde u$ and using the triangle inequality, we obtain  
	\[
	\|u\|_{C^{2+\alpha,1+\alpha/2}(\overline Q_T)}
	\le C(n,s,\alpha,\Omega,T)\bigl(\|f\|_{C^{\alpha,\alpha/2}(\overline Q_T)}+\|u_0\|_{C^{2+\alpha}(\overline\Omega)}+\|u\|_{L^\infty(Q_T)}\bigr).\qedhere
	\]  
\end{proof}

\subsection{Global-in-time existence}\label{sec:Global-in-time existence}
	By virtue of the blow-up criterion \cite[Theorem~2.2, Section~6.2]{Paz83}, Theorems~\ref{thm:existence} and~\ref{thm: reg up to bound} together guarantee that the classical solution of problem \eqref{problem}--\eqref{bdry-initial} exists globally in time. In other words, no finite-time blow-up occurs and the solution can be extended to any prescribed horizon.
	

{
	\begin{proof}[Proof of Corollary~\ref{cor:global-exist-reg}]
		The proof relies on a priori estimate \eqref{reg-up-to-bd} and the abstract theory for evolution equations. We proceed by steps.
		
		\textbf{Step 1: Local existence and regularity.} 
		Consider the operator $\mathcal{L}=-\Delta+(-\Delta)^s$ in $X:=C_0(\overline\Omega)$ with domain 
		\[
		D(A):=\left\{ ~u\in C^{2+\alpha}(\overline\Omega): u|_{\mathbb{R}^n \setminus \Omega}=0~\right\}.
		\]
		By the Hille--Yosida theorem (cf. \cite[Thm.~4.9]{Gri09}), the operator $-\mathcal{L}$ generates an analytic compact semigroup on $X$. 
		By Theorem~2.2 in \cite[Section~6.2]{Paz83}, for every $u_0\in C^{2+\alpha}(\overline\Omega)$ with $u_0|_{\mathbb{R}^n \setminus \Omega}=0$, there exists a unique mild solution $u$ defined on a maximal interval of existence $[0,t_{\max})$.
		
		By the regularity result Theorem~\ref{thm: reg up to bound}, this mild solution satisfies $u\in C^{2+\alpha,1+\alpha/2}(\overline{Q_{t}})$ for any $t<t_{\max}$, and the estimate \eqref{reg-up-to-bd} holds on each cylinder $Q_t$.
		
		\textbf{Step 2: A priori bound and global existence.}
		We claim that $t_{\max}=+\infty$. Suppose, by contradiction, that $t_{\max}< +\infty$. Then Theorem~2.2 of \cite[Section~6.2]{Paz83} implies
		\begin{equation}\label{blowup}
			\lim_{t\to t_{\max}}\|u(t)\|_{X}=+\infty.
		\end{equation}
		
		However, applying the a priori estimate \eqref{reg-up-to-bd} on $Q_{t_{\max}-\varepsilon}$ for small $\varepsilon>0$, we obtain
		\[
		\|u\|_{C^{2+\alpha,1+\alpha/2}(\overline{Q_{t_{\max}-\varepsilon}})}
		\leq C\left(\|f\|_{C^{\alpha,\alpha/2}(\overline Q_{t_{\max}})}+\|u_0\|_{C^{2+\alpha}(\overline\Omega)}+\|u\|_{L^\infty(Q_{t_{\max}-\varepsilon})}\right),
		\]
		where $C=C(n,s,\alpha,\Omega,t_{\max})$ is independent of $\varepsilon$. 
		
		Since $f\in C^{\alpha,\alpha/2}(\overline Q_\infty)$ and $u_0\in C^{2+\alpha}(\overline\Omega)$ are fixed, Theorem~\ref{Thm:global boundedness} yields
		\[
		\sup_{t\in[0,t_{\max})}\|u(t)\|_{L^\infty(Q_{t_{\max}-\varepsilon})} \leq C\bigl(\|f\|_{L^\infty(Q_{t_{\max}})}+\|u_0\|_{L^\infty(\Omega)}\bigr)< +\infty.
		\]
	    Consequently,
		\[
		\sup_{t\in[0,t_{\max})}\|u(t)\|_{X}< +\infty,
		\]
		which contradicts \eqref{blowup}. Therefore $t_{\max}=+\infty$.
		
		\textbf{Step 3: Regularity on the infinite cylinder.}
		For any $T>0$, the solution $u$ satisfies $u\in C^{2+\alpha,1+\alpha/2}(\overline Q_T)$ by Theorem~\ref{thm: reg up to bound}. The a priori estimate \eqref{reg-up-to-bd} is uniformly bounded on any compact time interval, so the regularity extends to all time \(t\in[0, +\infty)\).
	\end{proof}
	}

\section{Asymptotic Radial Symmetry and monotonicity}\label{sec: Asymptotic Radial Symmetry of Classical Solutions}
In this section, combining an adapted antisymmetric maximum principle, a narrow‑region lemma for $\partial_t+\mathcal{L}$, and the moving plane method, we will show the asymptotic radial symmetry of solutions to the problem \eqref{symmetry problem}.

We first introduce some notations that will be used throughout this section.
\begin{itemize}
	\item Let  
	\( T_\lambda = \{x \in \mathbb{R}^n \mid x_1 = \lambda, \text{ for } \lambda \in \mathbb{R}\} \) 
	be the moving planes.
	
	\item Define $\Sigma_\lambda = \{x \in \mathbb{R}^n \mid x_1 < \lambda\}$ as the region left to $T_\lambda$.
	
	\item For any point $x = (x_1,\ldots,x_n)$, let \( x^\lambda = (2\lambda - x_1, x_2, \cdots, x_n) \)
	be the reflection of \( x \) about the plane \( T_\lambda \).
	
	\item Denote
	\( u_\lambda(x,t) = u(x^\lambda, t) \) and \( w_\lambda(x,t) = u_\lambda(x,t) - u(x,t). \)
	Obviously, \( w_\lambda(x,t) \) is anti-symmetric with respect to the plane \( T_\lambda \).
	
	\item Define the \(\omega\)-limit set of \(u\) as
	\( \omega (u)= \{\, \varphi = \lim u(\cdot,t_k) \text{ for some } t_k \rightarrow +\infty \, \}\).
	
	\item For each \( \varphi(x) \in \omega(u) \), denote
	\( \psi_\lambda(x) = \varphi(x^\lambda) - \varphi(x) = \varphi_\lambda(x) - \varphi(x), \)
	which is an \( \omega \)-limit of \( w_\lambda(x, t) \).
\end{itemize}


In Subsections~\ref{sec:Asymptotic narrow region}-\ref{sec:Asymptotic maximum principles}, we present several key tools—asymptotic narrow region principle (Theorem \ref{thm: Asymptotic narrow region principle}), asymptotic maximum principle for antisymmetric functions (Lemma~\ref{lem: Asymptotic maximum principle for antisymmetric functions})  and asymptotic strong maximum principle for antisymmetric functions (Theorem \ref{thm: Asymptotic strong maximum principle for antisymmetric functions}). We then complete the proof of Theorem~\ref{thm:symmetry} in Subsection~\ref{sec:proof of theorem4}.

We point out that although only the case \(c_\lambda=0\) is needed in the proof of Theorem~\ref{thm:symmetry}, we state these results under the natural assumption that \(c_\lambda\) is merely bounded so that they remain directly applicable to nonlinear problems such as \(f(t,u)=u^p-u\) in our forthcoming work.

\subsection{Asymptotic narrow region principle}\label{sec:Asymptotic narrow region}
Let us first prove the following asymptotic narrow region principle.
\begin{theorem}[Asymptotic narrow region principle]\label{thm: Asymptotic narrow region principle}
	Let \(U\) be a bounded domain containing in the narrow slab
	\( \{x \mid \lambda - l < x_1 < \lambda\} \subset \Sigma_\lambda \)
	with some small \(l\).
	
	For \(\bar{t}\) sufficiently large, assume that \(w_\lambda(x,t) \in \big(C_{loc}^{2}(U) \cap L_{2s}^1\big) \times C^1([ \bar{t},+\infty))\) is uniformly bounded and lower semi-continuous in \(x\) on \(\overline{U}\), satisfying
	\begin{equation}\label{31}
		\begin{cases}
			\partial_t w_\lambda (x,t)  + \mathcal{L} w_\lambda(x,t) = c_\lambda(x,t) w_\lambda(x,t), & (x,t) \in U \times [\bar{t}, +\infty), \\
			w_\lambda(x,t) \geq 0, & (x,t) \in (\Sigma_\lambda \backslash U) \times [\bar{t}, +\infty), \\
			w_\lambda(x,t) = -w_\lambda(x^\lambda,t), & (x,t) \in U \times [\bar{t},+\infty),
		\end{cases}
	\end{equation}
	where \(c_\lambda(x,t)\) is bounded from above. Then, for sufficiently small \(l\), we have
	\begin{equation*}
		\liminf_{t \to +\infty} w_\lambda(x,t) \geq 0, \quad \forall \, x \in U.
	\end{equation*}
\end{theorem}	

\begin{proof}
	Fix \( m >0\) to be chosen later and set  
	\( \tilde{w}_\lambda(x,t) = e^{mt} w_\lambda(x,t) \), then
	\begin{equation}\label{32}
		\partial_t \tilde{w}_\lambda(x,t) + (-\Delta)^s \tilde{w}_\lambda(x,t) = (m + c_\lambda(x,t)) \tilde{w}_\lambda(x,t), \quad (x,t) \in U \times [\bar{t}, +\infty). 
	\end{equation}
	
	We first claim that for any \( T > \bar{t}, \)  
	\begin{equation}\label{33}
		\tilde{w}_\lambda(x,t) \geq \min\{0, \inf_{U} \tilde{w}_\lambda(x,\bar{t})\}, \quad (x,t) \in U \times [\bar{t}, T].
	\end{equation}
	
	We proceed by contradiction. 
	If \eqref{33} is false, by \eqref{31} and the lower semi-continuity of \( w_\lambda \), there exists \((x_0, t_0) \in U \times (\bar{t}, T]\) such that  
	\begin{equation*}
		\tilde{w}_\lambda(x_0, t_0) = \min_{\Sigma_\lambda \times (\bar{t}, T]} \tilde{w}_\lambda(x,t) < \min\{0, \inf_{U} \tilde{w}_\lambda(x,\bar{t})\}.
	\end{equation*}
	
	Then $ -\Delta \tilde{w}_\lambda(x_0, t_0)\leq 0$ and 
	\begin{equation}\label{35}
		\partial_t \tilde{w}_\lambda(x_0, t_0) \leq 0 .
	\end{equation} 
	Thus, we obtain
	\begin{equation}\label{37}
		\begin{aligned}
			&\quad \left(-\Delta + (-\Delta)^s \right)\tilde{w}_\lambda(x_0, t_0)\leq (-\Delta)^s \tilde{w}_\lambda(x_0, t_0) \\
			&= c_{n,s} P.V. \int_{\Sigma_\lambda} \frac{\tilde{w}_\lambda(x_0, t_0) - \tilde{w}_\lambda(y, t_0)}{|x_0 - y|^{n+2s}} dy  + c_{n,s} \int_{\mathbb{R}^n \setminus \Sigma_\lambda} \frac{\tilde{w}_\lambda(x_0, t_0) - \tilde{w}_\lambda(y, t_0)}{|x_0 - y|^{n+2s}} dy \\
			&= c_{n,s} P.V. \int_{\Sigma_\lambda} \frac{\tilde{w}_\lambda(x_0, t_0) - \tilde{w}_\lambda(y, t_0)}{|x_0 - y|^{n+2s}} dy  + c_{n,s} \int_{\Sigma_\lambda} \frac{\tilde{w}_\lambda(x_0, t_0) - \tilde{w}_\lambda(y^\lambda, t_0)}{|x_0 - y^\lambda|^{n+2s}} dy \\
			&\leq c_{n,s} \int_{\Sigma_\lambda} \frac{2\tilde{w}_\lambda(x_0, t_0)}{|x_0 - y^\lambda|^{n+2s}} dy.
		\end{aligned}
	\end{equation}
	Selecting a subregion $H \subset \Sigma_\lambda$ defined by
	\[ H = \left\{ ~y = (y_1, y') \in \Sigma_\lambda ~:~ l < y_1 - (x_0)_1 < 2l, \, |y' - x_0'| < l ~ \right\}, 
	\]
	then we have
	\begin{equation}\label{38}
		\int_{\Sigma_\lambda} \frac{1}{|x_0 - y^\lambda|^{n+2s}} dy \geq \int_H \frac{1}{|x_0 - y^\lambda|^{n+2s}} dy \geq \frac{C}{l^{n+2s}} |H| = \frac{C}{l^{2s}}.
	\end{equation}
	
	Combining \eqref{32} and \eqref{37}-\eqref{38}, we obtain
	\begin{equation}\label{36}
		\begin{aligned}
			&\quad \partial_t \tilde{w}_\lambda(x_0, t_0) = -(-\Delta+(-\Delta)^s) \tilde{w}_\lambda(x_0, t_0) + (m + c_\lambda(x_0, t_0)) \tilde{w}_\lambda(x_0, t_0)\\
			&\geq -c_{n,s} \int_{\Sigma_\lambda} \frac{2\tilde{w}_\lambda(x_0, t_0)}{|x_0 - y^\lambda|^{n+2s}} dy + (m + c_\lambda(x_0, t_0)) \tilde{w}_\lambda(x_0, t_0)\\
			&\geq \left( -C/l^{2s} + m + c_\lambda(x_0, t_0) \right) \tilde{w}_\lambda(x_0, t_0),
		\end{aligned}
	\end{equation}
	due to $\tilde{w}_\lambda(x_0, t_0) < 0$.
	
	Since \(c_\lambda(x, t)\) is bounded from above, we have $c_\lambda(x_0, t_0) < \frac{C}{2l^{2s}}$ for sufficiently small $l$.
	Taking $m = \frac{C}{2l^{2s}} $, we conclude that the right-hand side of \eqref{36} is strictly positive, which contradicts \eqref{35}. 
	
	Thus, we have proven the claim \eqref{33}.
	
	Therefore, there exists \( C_1 > 0 \) such that
	\begin{equation*}
		\tilde{w}_{\lambda}(x,t) \geq \min\left\{0, \inf_{U} \tilde{w}_{\lambda}(x,\bar{t}) \right\} \geq -C_1, \quad  (x,t) \in U \times [\bar{t},T], \, \forall~ T > \bar{t}. 
	\end{equation*}
	Consequently,
	\[
	w_{\lambda}(x,t) \geq -C_1 e^{-mt}, \quad \forall ~t > \bar{t}.
	\] 
	Taking the limit as \( t \to +\infty \), the theorem is proved. 
\end{proof}

\subsection{Asymptotic maximum principles}\label{sec:Asymptotic maximum principles}
We now ready to establish both asymptotic maximum principles for antisymmetric functions below.
\begin{lemma}[Asymptotic maximum principle for antisymmetric functions]\label{lem: Asymptotic maximum principle for antisymmetric functions}
	Let $U$ be a bounded domain in $\Sigma_\lambda$. Assume that $v(x,t) \in (C^{2}_{loc}(U) \cap L_{2s}^1) \times C^1([0, +\infty))$ is lower semi-continuous in $x$ on $\bar{U}$ and satisfies
	\[
	\begin{cases}
		\partial_t v + \mathcal{L} v \geq c_\lambda(x,t) v, & (x,t) \in U \times [0, +\infty), \\
		v(x^\lambda,t) = -v(x,t), & (x,t) \in \Sigma_\lambda \times [0, +\infty), \\
		v(x,t) \geq 0, & (x,t) \in (\Sigma_\lambda \backslash U) \times [0, +\infty), \\
		v(x,0) \geq 0, & x \in U.
	\end{cases}
	\]
	If $c_\lambda(x,t)$ is bounded from above, then
	\begin{equation}\label{39}
		v(x,t) \geq 0, \quad (x,t) \in U \times [0,T], \, \forall \, T > 0.
	\end{equation}
\end{lemma}
\begin{proof}
	Since $c_{\lambda }(x,t)$ is bounded from above, we can choose $m < 0$ such that $m + c_{\lambda}(x,t) < 0$. Let
	$\tilde{v}(x,t) = e^{mt}v(x,t)$, then $\tilde{v}(x,t)$ satisfies
	\begin{equation}\label{34}
		\partial_t \tilde{v} -\Delta \tilde{v}+ (-\Delta)^s \tilde{v} \geq (m + c_{\lambda}(x,t))\tilde{v}, \, (x,t) \in U \times [0, +\infty).
	\end{equation}
	
	To prove \eqref{39}, it suffices to prove 
	\begin{equation}\label{30}
		\tilde{v}(x,t) \geq \inf_{U} \tilde{v}(x,0) \geq 0 \quad \text{in } \, U \times [0,T].
	\end{equation}
	We proceed by contradiction. 
	If \eqref{30} is false, there exists $(x_0, t_0) \in U \times (0,T]$ such that
	\[
	\tilde{v}(x_0, t_0) = \min_{\mathbb{R}^n \times [0,T]} \tilde{v}(x,t) < 0,
	\]
	then
	$(m + c_\lambda(x_0, t_0))\tilde{v}(x_0, t_0) > 0$.
	Since $(x_0,t_0)$ is the minimum point, we obtain $\partial_t \tilde{v}(x_0, t_0) \leq 0$ and
	\begin{equation*}
		\begin{aligned}
			&\quad (-\Delta+(-\Delta)^s) \tilde{v}(x_0, t_0) \leq c_{n,s} P.V. \int_{\mathbb{R}^n} \frac{\tilde{v}(x_0,t_0)-\tilde{v}(y,t_0)}{|x_0-y|^{n+2s}}dy\\
			&=c_{n,s} P.V. \left(\int_{\Sigma_\lambda} \frac{\tilde{v}(x_0,t_0)-\tilde{v}(y,t_0)}{|x_0-y|^{n+2s}}dy+\int_{\mathbb{R}^n \setminus \Sigma_\lambda} \frac{\tilde{v}(x_0,t_0)-\tilde{v}(y,t_0)}{|x_0-y|^{n+2s}}dy\right)\\
			&\leq c_{n,s} P.V. \left( \int_{\Sigma_\lambda} \frac{\tilde{v}(x_0,t_0)-\tilde{v}(y,t_0)+\tilde{v}(x_0,t_0)-\tilde{v}(y^\lambda,t_0)}{|x_0-y|^{n+2s}} \right) dy\\
			&= c_{n,s} P.V. \int_{\Sigma_\lambda} \frac{2\tilde{v}(x_0,t_0)}{|x_0-y|^{n+2s}} dy <0.
		\end{aligned}
	\end{equation*}
	This contradicts \eqref{34}. Thus, \eqref{30} holds.
\end{proof}

\begin{theorem}[Asymptotic strong maximum principle for antisymmetric functions]\label{thm: Asymptotic strong maximum principle for antisymmetric functions}
	Let \(s\in(0,1/2)\). Denote \(\Omega_\lambda:=\Sigma_\lambda \cap B_1(0)\).
	For sufficiently large $\bar{t}$, suppose that $w_\lambda(x,t) \in C^{2,1}(\Omega_\lambda \times[\bar{t}, +\infty)) $ is bounded, vanishes outside \(B_1(0)\), and satisfies
	\begin{equation}\label{41}
		\begin{cases}
			\partial_t w_\lambda + \mathcal{L} w_\lambda = 0, & (x,t) \in \Omega_\lambda \times [\bar{t}, + \infty), \\
			w_\lambda(x,t) = -w_\lambda(x^\lambda, t), & (x,t) \in \Omega_\lambda \times [\bar{t}, + \infty), \\
			\liminf\limits_{t \to +\infty} w_\lambda(x,t) \geq 0, & x \in \Omega_\lambda,
		\end{cases}
	\end{equation}
	with $c_\lambda(x,t)$ is bounded. 
	
	Assume $\psi_\lambda > 0$ somewhere in $\Omega_\lambda$. Then $\psi_\lambda(x) > 0$ in $\Omega_\lambda$.
\end{theorem}

\begin{proof}
	The proof follows via an asymptotic approach as in \cite[Theorem 2.7]{CWNH21}. For the reader's convenience, we nevertheless provide the complete proof.
	
	For any $\psi_\lambda \in \omega(w_\lambda)$, there exists $t_k$ such that $w_\lambda(x, t_k) \to \psi_\lambda(x)$ as $t_k \to +\infty$. Let 
	\[ w_k(x, t) = w_\lambda(x, t + t_k - 1).\]
	Then
	\[
	\partial_t w_k(x, t) + \mathcal{L} w_k(x, t) = 0, \quad (x, t) \in \Omega_\lambda \times [\bar{t}, +\infty).
	\] 
	
	By Theorem~\ref{thm: reg up to bound}, we infer that there is a subsequence of $w_k(x, t)$ (still denoted by $w_k(x, t)$) which converges uniformly to a function $w_\infty(x, t)$ in $\Omega_\lambda \times [0, 2]$ and
	\begin{equation}
		\begin{aligned}
			\frac{\partial w_k}{\partial t}(x,t)  + \mathcal{L} w_k(x,t) \to \frac{\partial w_\infty}{\partial t}(x,t) + \mathcal{L} w_\infty(x,t), \quad \text{as } k \to +\infty.
		\end{aligned}
	\end{equation} 
	Particularly, we get
	\[
	w_\lambda(x,t_k) = w_k(x,1) \to w_\infty(x,1) = \psi_\lambda(x) \quad \text{as } k \to +\infty.
	\]
	
	Take $m > 0$. Consider the function $ \tilde{w}(x,t) = e^{mt} w_\infty(x,t)$.
	By the third condition in \eqref{41}, we have
	$	\tilde{w}(x,t) \geq 0$ in $ \Omega_\lambda \times [0,2] $
	and thus
	\begin{equation}\label{46}
		\partial_t \tilde{w}(x,t) + (-\Delta + (-\Delta)^s) \tilde{w}(x,t) = m\tilde{w}(x,t) \geq 0, \quad \text{in } \Omega_\lambda \times [0,2].
	\end{equation}
	
	Since $\psi_\lambda > 0$ somewhere in $\Omega_\lambda$, by continuity, there exists a set $D \Subset  \Omega_\lambda$ such that
	\begin{equation}\label{47}
		\psi_\lambda(x) > C > 0, \quad x \in D.
	\end{equation}
	We next prove $\psi_\lambda(x) > 0$ in $\Omega_\lambda$ by distinguishing the two cases $x \in D$ and $x \in \Omega_\lambda \setminus D$.
	
	\textit{Case 1.} For $x\in D$, by the continuity of $w_\infty$, there exists $0 < \varepsilon_0 < 1$, such that
	\begin{equation}\label{42}
		w_\infty(x,t) > C/2, \quad (x,t) \in D \times [1 - \varepsilon_0, 1 + \varepsilon_0].
	\end{equation}
	
	\textit{Case 2.} For $\bar{x} \in \Omega_\lambda \backslash D$, choose $\delta = \min\{\text{dist}(\bar{x}, D), \text{dist}(\bar{x},\Omega_\lambda)\} > 0$, then $B_\delta(\bar{x}) \subset \Omega_\lambda \backslash D$.
	
	We employ the auxiliary function 
	\[
	\underline{w}(x,t) = \chi_{D \cup D_\lambda}(x)\tilde{w}(x,t) + \epsilon \eta(t)g(x),
	\]
	introduced in \cite[Theorem 2.7]{CWNH21}. Here $D_\lambda$ is the reflection of $D$ about the plane $T_\lambda$, 
	$\eta(t)\in C^\infty(\mathbb R)$ is a cutoff function satisfying
	\[
	\eta(t)=
	\begin{cases}
		1, \quad t\in [1 - \varepsilon_0, 1 + \varepsilon_0/2],\\
		0, \quad t\not \in [1 - \varepsilon_0, 1 + \varepsilon_0],
	\end{cases}
    \]
	and 	
	$g(x) = (\delta^2 - |x - \bar{x}|^2)_+^s - (\delta^2 - |x - \bar{x}^\lambda|^2)_+^s$. 
	As in \cite{CWNH21}, one can deduce that for each fixed $t \in [1 - \varepsilon_0, 1 + \varepsilon_0]$ and for any $x \in B_\delta(\bar{x})$
	\[ (-\Delta)^s g(x) \leq C_0 \quad \text{and} \quad (-\Delta)^s(\chi_{D \cup D_\lambda}(x) \tilde{w}(x, t))\leq -C_1.\]
	Moreover, for any $x \in B_\delta(\bar{x})$, \(-\Delta g(x)\le C_2\).
	
	Thus, we derive for any $(x, t) \in B_\delta(\bar{x}) \times [1 - \varepsilon_0, 1 + \varepsilon_0]$
	\begin{equation}\label{45}
		\begin{aligned}
			&\quad \partial_t \underline{w}(x, t) + \mathcal{L} \underline{w}(x, t)\\
			&=\chi_{D \cup D_\lambda}(x) \partial_t \tilde{w}(x,t) + \epsilon \eta'(t)g(x) + \mathcal{L} (\chi_{D \cup D_\lambda}(x)\tilde{w}(x,t)) + \epsilon \eta(t) \mathcal{L}  g(x)\\
			& = \varepsilon \eta'(t) g(x) + (-\Delta)^s (\chi_{D \cup D_\lambda} (x)\tilde{w}(x, t)) + \varepsilon \eta(t)(-\Delta) g(x) + \varepsilon \eta(t)(-\Delta)^s g(x)\\
			&\leq \varepsilon \eta'(t) g(x) - C_1 + \varepsilon \eta(t) (C_0 + C_2).
		\end{aligned}
	\end{equation}
	Choose $\varepsilon > 0$ sufficiently small such that 
	\[
	\partial_t \underline{w}(x,t) + (-\Delta+(-\Delta)^s) \underline{w}(x,t) \leq 0, \quad (x,t) \in B_\delta(\bar{x}) \times [1 - \varepsilon_0, 1 + \varepsilon_0].
	\]
	
	Set $v(x,t) = \tilde{w}(x,t) - \underline{w}(x,t)$. Then $v(x,t) = -v(x^\lambda,t)$ and by \eqref{46} and \eqref{45} 
	\[
	\partial_t v(x,t) + \mathcal{L} v(x,t) \geq 0, \quad (x,t) \in B_\delta(\bar{x}) \times [1 - \varepsilon_0, 1 + \varepsilon_0].
	\]
	Also, by the definition of $\underline{w}(x,t)$, we have
	\[
	v(x,t) \geq 0 \text{ in } (\Omega_\lambda \setminus B_\delta(\bar{x})) \times [1 - \varepsilon_0, 1 + \varepsilon_0] \, \text{ and } \, v(x,0) \geq 0 \text{ in } B_\delta(\bar{x}).
	\]
	
	Applying Lemma \ref{lem: Asymptotic maximum principle for antisymmetric functions}, we obtain
	\[
	v(x,t) \geq 0, \quad (x,t) \in B_\delta(\bar{x}) \times [1 - \varepsilon_0, 1 + \varepsilon_0].
	\]
	In particular, taking $t=1$ we have
	\(
	w_{\infty}(x,1) \geq e^{-m} \varepsilon g(x)
	\)
	for any \(x \in B_\delta(\bar{x})\).
	Since $g(\bar{x}) = \delta^{2s}$, we conclude
	\begin{equation}\label{48}
		\psi_\lambda(\bar{x}) = w_{\infty}(\bar{x},1) \geq e^{-m} \varepsilon \delta^{2s} > 0. 
	\end{equation}
	
	Since $\bar x\in \Omega_\lambda \setminus D$ is arbitrary, \eqref{47} and \eqref{48} yield $\psi_\lambda>0$ in $\Omega_\lambda$.
\end{proof}

\subsection{Proof of Theorem~\ref{thm:symmetry}}\label{sec:proof of theorem4}
	We remark that the assumption of the existence of a positive solution $u\in  C^{2+\alpha,1+\alpha/2}(\overline Q_\infty)$ is justified by Corollary~\ref{cor:global-exist-reg} and Theorem~\ref{thm:max_principle}.
	
\begin{proof}[Proof of Theorem~\ref{thm:symmetry}]
	For any $\varphi \in \omega(u)$, if $\varphi(x) \equiv 0$, the conclusion is trivial. Without loss of generality, we assume  $\varphi \not \equiv 0$. 
	Fix a direction, say $x_1$, and define for $\lambda\in(-1,0)$, \(\Omega_\lambda=\Sigma_\lambda\cap B_1(0)\) and $w_\lambda$ satisfies
	\begin{equation}\label{49}
		\begin{cases}
			\partial_t w_\lambda + (-\Delta+(-\Delta)^s) w_\lambda = 0, & (x,t) \in \Omega_\lambda \times (0,+\infty), \\
			w_\lambda(x,t) = -w_\lambda(x^\lambda,t), & (x,t) \in \Omega_\lambda \times (0, +\infty).
		\end{cases}
	\end{equation}
	We carry out the moving-plane procedure in three steps.
	
	\medskip
	\noindent\textit{Step 1: Starting.} For $\lambda$ close to $-1$, $\Omega_\lambda$ is a narrow region. Since $u(x,t) = 0$ in $ B_1^c(0) \times (0, +\infty)$, 
	\(
	w_\lambda(x,t) \geq 0, \, (x,t) \in (\Sigma_\lambda \backslash \Omega_\lambda)\times (0, +\infty).
	\)
	Theorem~\ref{thm: Asymptotic narrow region principle} yields
	\[
	\psi_\lambda(x)\ge 0\quad\text{for all }x\in\Omega_\lambda\text{ and all }\varphi\in\omega(u).
	\]
	
	\medskip
	\noindent\textit{Step 2: Moving.}  
	Define
	\[
	\lambda_0=\sup\bigl\{\lambda\le 0:\psi_\mu\ge 0\text{ on }\Omega_\mu\text{ for all }\mu\le\lambda\text{ and all }\varphi\in\omega(u)\bigr\}.
	\]
	We will show $\lambda_0=0$.  
	Assume for contradiction that $\lambda_0<0$.  
	By the definition of $\lambda_0$, $\psi_{\lambda_0}(x)\ge 0$ on $\Omega_{\lambda_0}$ for every $\varphi\in\omega(u)$.  	
	
	We first claim that for any $\varphi \in \omega(u)$, there exists $x_\varphi \in \Sigma_{\lambda_0}$ such that
	\begin{equation}\label{51}
		\psi_{\lambda_0}(x_\varphi) > 0. 
	\end{equation}
	
	If \eqref{51} is false, there exists some $\bar{\varphi} \in \omega(u)$ such that
	\(
	\bar{\psi}_{\lambda_0}(x) \equiv 0 
	\)
	in \(\Sigma_{\lambda_0}\). By the outer condition of $u$, we have $\bar{\varphi}(x) \equiv 0$ in $B_1^c(0) \cap \Sigma_{\lambda_0}$, therefore, there exists $x_0 \in B_1(0)$ such that 
	\begin{equation}\label{52}
		\bar{\varphi}(x_0)= \bar{\varphi}_{\lambda_0}({x_0}^\lambda) = 0.
	\end{equation}
	For this $\bar{\varphi}$, there exists $t_k$ such that $u(x,t_k) \to \bar{\varphi}(x)$ as $t_k \to +\infty$. Set
	\[u_k(x,t):= u(x,t+t_k-1). \]
	Then $u_k(x,1)=u(x,t_k) \to \bar{\varphi}(x)$ as $k \to +\infty$ and $\partial_t u_k +(-\Delta+(-\Delta)^s)u_k =0$.
	
	By the regularity estimates in Theorem~\ref{thm: reg up to bound}, there exists some function $u_\infty(x, t)$ such that $u_k \to u_\infty \geq 0$ and 
	\begin{equation}\label{53}
		\partial_t u_\infty (x,t)  + (-\Delta + (-\Delta)^s) u_\infty(x,t) =0, \quad (x,t)\in \Omega_\lambda \times[0,2].
	\end{equation} 
	Moreover, by \eqref{52},  $u_\infty(x_0,1) = \bar{\varphi}(x_0)=0$. Thus $\partial_t u_\infty(x_0,1) = 0$, $-\Delta u_\infty(x_0,1) \leq 0$ and
	\[
	(-\Delta)^s u_\infty(x_0,1) = c_{n,s} P.V. \int_{B_1(0)} \frac{-u_\infty(y,1)}{|x_0-y|^{n+2s}} dy < 0,
	\]
	where we used in the last inequality that $u_\infty(y,1) \neq 0$ in $B_1(0)$ (since $\bar{\varphi} \neq 0$ in $B_1(0)$). This contradicts \eqref{53}. Thus, we obtain the claim \eqref{51}.
	
	Using the asymptotic strong maximum principle for antisymmetric functions (Theorem~\ref{thm: Asymptotic strong maximum principle for antisymmetric functions}), we have $
	\psi_{\lambda_0}(x) > 0 $ on \(\Omega_{\lambda_0}\) for all $\varphi \in \omega(u)$.
	
	Thus for any $\delta > 0$ small, for each $\psi_{\lambda_0}=\varphi(x^\lambda)-\varphi(x)$, there exists a constant $C_\varphi > 0$ (depending on $\varphi$), such that \(\psi_{\lambda_0}(x) \geq C_\varphi\) on \(\overline{\Omega_{\lambda_0-\delta}}\).
	By the compactness of the $\omega$-limit set, we claim further that, for all $\varphi \in \omega(u)$, there exists a universal constant $C_0$ such that 
	\begin{equation}\label{54}
		\psi_{\lambda_0}(x) \geq C_0>0, \, x\in \overline{\Omega_{\lambda_0-\delta}}.
	\end{equation}
	By the continuity of $\psi_\lambda$ with respect to $\lambda$, for each $\psi_\lambda$ (corresponding to $\varphi \in \omega(u)$), there exists $\varepsilon_\varphi >0$ (depending on $\varphi$) such that
	\[
	\psi_\lambda(x) \geq C_0/2 > 0, \quad x \in \overline{\Omega_{\lambda_0-\delta}}, \quad \forall \; \lambda \in (\lambda_0, \lambda_0 + \varepsilon_\varphi).
	\]
	
	Through a similar compactness argument as in deriving \eqref{54}, we conclude that, for all $\psi_\lambda$, there exists a universal $\varepsilon > 0$ such that
	\begin{equation}\label{55}
		\psi_\lambda(x) \geq C_0/2 > 0, \quad x \in \overline{\Omega_{\lambda_0-\delta}}, \quad \forall \; \lambda \in (\lambda_0, \lambda_0 + \varepsilon). 
	\end{equation}
	As a consequence, for $t$ sufficiently large, we have
	\[
	w_\lambda(x,t) \geq 0, \quad x \in \overline{\Omega_{\lambda_0-\delta}}, \quad \forall \; \lambda \in (\lambda_0, \lambda_0 + \varepsilon).
	\]
	
	Since $\delta > 0$ is small, by \eqref{55}, we can choose $\varepsilon > 0$ small, such that $\Omega_\lambda \backslash \Omega_{\lambda_0-\delta}$ is a narrow region for $\lambda \in (\lambda_0, \lambda_0 + \varepsilon)$. Applying the asymptotic narrow region principle (Theorem \ref{thm: Asymptotic narrow region principle}), we arrive at
	\begin{equation}\label{56}
		\psi_{\lambda}(x) \geq 0, \quad \forall x \in \Omega_{\lambda} \setminus \Omega_{\lambda_0 - \delta}.
	\end{equation}
	
	Combining \eqref{55} and \eqref{56}, we conclude that
	\[
	\psi_{\lambda}(x) \geq 0, \quad \forall x \in \Omega_{\lambda}, \quad \forall \lambda \in (\lambda_0, \lambda_0 + \epsilon), \quad \forall \varphi \in \omega(u).
	\]
	This contradicts the definition of $\lambda_0$. Therefore, $\lambda_0 = 0$. 
		
	\medskip
	\noindent\textit{Step 3: Radial monotonicity.} It follows that for all $\varphi \in \omega(u)$ and $x \in \Omega_0$, $\psi_0(x) \geq 0$.
	Equivalently, 
	\begin{equation}\label{57}
		\varphi(-x_1, \cdots, x_N) \leq \varphi(x_1, \cdots, x_N), \ \text{for all } \varphi \in \omega(u),\ 0 < x_1 < 1. 
	\end{equation}
	Since the direction $x_1$ is arbitrary,  \eqref{57} implies that all $\varphi\in\omega(u)$ is radially symmetric.  
	Strict decrease follows from the strong positivity of $\psi_\lambda$ for $-1<\lambda<0$.
\end{proof}

\begin{description}	
\item [\textbf{Conflict of interest}]
The authors have no conflict of interest to disclose.

\item [\textbf{Data availability}] 
No data was used for the research described in the article.
\end{description}

\bibliographystyle{alpha}

\bibliography{parabolic}
	
\end{document}